\documentclass[11pt]{article}
\usepackage{amssymb,amsmath,latexsym}
\begin{document}

\begin{doublespace}

\newtheorem{thm}{Theorem}[section]
\newtheorem{lemma}[thm]{Lemma}
\newtheorem{defn}[thm]{Definition}
\newtheorem{prop}[thm]{Proposition}
\newtheorem{corollary}[thm]{Corollary}
\newtheorem{cor}[thm]{Corollary}
\newtheorem{remark}[thm]{Remark}
\newtheorem{example}[thm]{Example}
\numberwithin{equation}{section}
\def\ee{\varepsilon}
\def\qed{{\hfill $\Box$ \bigskip}}
\def\NN{{\cal N}}
\def\AA{{\cal A}}
\def\MM{{\cal M}}
\def\BB{{\cal B}}
\def\LL{{\cal L}}
\def\DD{{\cal D}}
\def\FF{{\cal F}}
\def\EE{{\cal E}}
\def\EE{{\cal E}}
\def\QQ{{\cal Q}}
\def\RR{{\mathbb R}}
\def\R{{\mathbb R}}
\def\L{{\bf L}}
\def\E{{\mathbb E}}
\def\F{{\bf F}}
\def\C{{\bf C}}
\def\P{{\mathbb P}}
\def\N{{\mathbb N}}
\def\eps{\varepsilon}
\def\wh{\widehat}
\def\pf{\noindent{\bf Proof.} }

\title{\Large \bf  Boundary Behavior of Harmonic Functions for
Truncated Stable Processes}
\author{Panki Kim\thanks{The research of this author is supported by Research Settlement Fund for the new faculty of Seoul
National University.}\\
Department of Mathematics\\
Seoul National University\\
Seoul 151-742, Republic of Korea\\
Email: pkim@snu.ac.kr \smallskip \\
and
\smallskip \\
Renming Song\thanks{The research of this author is supported in part
by a joint
US-Croatia grant INT 0302167.}\\
Department of Mathematics\\
University of Illinois \\
Urbana, IL 61801, USA\\
Email: rsong@math.uiuc.edu }
\date{ }
\maketitle

\begin{abstract}
For any $\alpha\in (0, 2)$, a truncated symmetric $\alpha$-stable
process in $\R^d$ is a symmetric L\'evy process in $\R^d$ with no
diffusion part and with a  L\'evy density given by
$c|x|^{-d-\alpha}\, 1_{\{|x|< 1\}}$ for some constant $c$. In
\cite{KS} we have studied the potential theory of truncated
symmetric stable processes. Among other things, we proved that the
boundary Harnack principle is valid for the positive harmonic
functions of this process in any bounded
convex domain and showed that the Martin boundary of any bounded
convex domain with respect to this process is the same as the
Euclidean boundary. However, for truncated symmetric stable
processes, the boundary Harnack principle is not valid in non-convex
domains. In this paper, we show that, for a large class of not
necessarily convex bounded open sets in $\R^d$ called bounded
roughly connected $\kappa$-fat open sets (including bounded
non-convex $\kappa$-fat domains), the Martin boundary with respect
to any truncated symmetric stable process is still the same as the
Euclidean boundary. We also show that, for truncated symmetric
stable processes a relative Fatou type theorem is true in bounded
roughly connected $\kappa$-fat open sets.
\end{abstract}



\section{Introduction}

Recently there has been a lot of interest in studying the potential
theory of discontinuous stable processes due to their importance in
theory as well as applications. Many deep results have been
established. However in a lot of applications one needs to use
discontinuous Markov processes which are not stable processes. For
example, in mathematical finance, it has been observed that even
though discontinuous stable processes provide better representations
of financial data than Gaussian processes (cf. \cite{HPR}),
financial data tend to become more Gaussian over a longer time-scale
(see \cite{M} and the references therein). The so called
relativistic stable processes (see \cite{CS, Ry}) have this required
property: they behave like  discontinuous stable processes in small
scale and behave like Brownian motion in large scale. Other
processes having this kind of property can be obtained by
``tempering'' stable processes, that is, by multiplying the L\'evy
densities of stable processes with strictly positive and completely
monotone decreasing factors (see \cite{R}).

In \cite{KS}, we considered an extreme case of ``tempering'': we
truncated the L\'evy densities of stable processes and obtained a
class of L\'evy processes called truncated stable processes. For any
$\alpha\in (0, 2)$, a truncated symmetric $\alpha$-stable process is
a symmetric L\'evy process in $\R^d$ with no diffusion part and with
a L\'evy density $l(x)$ coincides with the L\'evy density of a
symmetric $\alpha$-stable process for $|x|$ small (say, $|x|<1$) and
is equal to zero for $|x|$ large (say, $|x|\ge 1$). In other words,
a truncated symmetric $\alpha$-stable process is a symmetric L\'evy
process in $\R^d$ with no diffusion part and with a  L\'evy density
given by $c|x|^{-d-\alpha}\, 1_{\{|x|< 1\}}$ for some positive
constant $c$. Truncated stable processes are very natural and
important in applications where only jumps up to a certain size are
allowed. In \cite{KS} we studied the potential theory of truncated
symmetric stable processes. Among other things, we proved that the
boundary Harnack principle is valid for the positive harmonic
functions of this process in bounded convex
domains and showed that the Martin boundary of any bounded convex
domain with respect to this process is the same as the Euclidean
boundary. However, for truncated symmetric stable processes, the
boundary Harnack principle is not valid in non-convex domains (see
the last section of \cite{KS} for a counterexample). A very natural
and very important question is: can one identify the Martin boundary
of bounded non-convex domains with respect to truncated symmetric
stable processes?

Recently, a relative Fatou type theorem has been established for
symmetric stable processes.
It is known that if $u$ and $h$ are positive harmonic function for
a symmetric $\alpha$-stable process in a bounded $\kappa$-fat open set $D$
with $h$ vanishing on $D^c$,
then the non-tangential limit of $u/h$ exists almost everywhere with
respect to the Martin measure of $h$. The assumption that $h$
vanishes on $D^c$ is necessary (see \cite{BY}). The above relative
Fatou type theorem was proved in \cite{BD} for bounded $C^{1,1}$
domains and extended to more general open sets in \cite{K2} and
\cite{MR} independently (see \cite{K} for a Fatou type theorem for
another class of discontinuous processes). With the recent results
obtained in \cite{KS} in hand, one naturally comes to the following
question: can one prove a relative Fatou type theorem for truncated
symmetric stable processes?

In this paper we will continue our study of truncated symmetric
stable processes. We will show that, for a large class of open sets
called bounded roughly connected $\kappa$-fat open sets (see
Definitions \ref{d:rc} and \ref{fat}), the Martin boundary with
respect to any truncated symmetric stable process is identical to
the Euclidean boundary and a relative Fatou type theorem holds. The
main tool for establishing this is the fact that, for any bounded
roughly connected $\kappa$-fat open set, the Green function of a
truncated symmetric $\alpha$-stable process is comparable to that of
a symmetric $\alpha$-stable process. This result on Green function
comparison is obtained by combining ideas from \cite{GR, KL, KS6}.

This paper is organized as follows. In section 2, we recall the
definition of a truncated stable process and collect some  basic
results on this process. Section 3 contains the result on the
comparison of Green functions which is used in Sections 4 to study
the Martin boundary. In Section 5, we establish a relative Fatou
type theorem for truncated stable processes in bounded roughly
connected $\kappa$-fat open sets. The main idea of Section 5 is
similar to that of \cite{K2}, which is inspired by Doob's approach.

Throughout this paper, for two real numbers $a$ and $b$, we denote
$a \wedge b := \min \{ a, b\}$ and $a \vee b := \max \{ a, b\}$. The
distance between $x$ and $\partial D$ is denote by $\rho_D (x)$. In
this paper, we use ``$:=$" to denote a definition, which is read as
``is defined to be". In this paper, the values and labeling of the
 constants $c, c_1, c_2, \cdots$ start anew in the statement of each
result.

\section{Truncated Stable Processes}

In this section we recall the definition of a truncated stable
process and collect some basic properties of this process from
\cite{KS}.

Throughout this paper we assume  $\alpha\in (0, 2)$ and $d\ge 2$.
Recall that a symmetric $\alpha$-stable process $X=(X_t, \P_x)$ in $\R^d$
is a L\'evy process such that
$$
\E_x\left[e^{i\xi\cdot(X_t-X_0)}\right]=e^{-t|\xi|^{\alpha}},
\quad \quad \mbox{ for every } x\in \R^d \mbox{ and } \xi\in \R^d.
$$
It is well known that
$$
|\xi|^{\alpha}=\int_{\R^d}(1-e^{i\xi \cdot x}+i\xi \cdot
x1_{\{|x|<1\}}) \nu^X(x)dx, \qquad  \xi\in \R^d,
$$
where
$$
\nu^X(x):={\cal A} (d, - \alpha)|x|^{-(d+\alpha)}, \qquad x\in \R^d,
$$
with ${\cal A}(d, -\alpha):= \alpha2^{\alpha}\pi^{-d/2}
\Gamma(\frac{d+\alpha}2) \Gamma(1-\frac{\alpha}2)^{-1}$. Here
$\Gamma$ is the Gamma function defined by $\Gamma(\lambda):=
\int^{\infty}_0 t^{\lambda-1} e^{-t}dt$ for every $\lambda > 0$.
$\nu^X$ is called the L\'evy density of $X$.

By a truncated symmetric $\alpha$-stable process in $\R^d$ we mean a
symmetric L\'evy process $Y=(Y_t, \P_x)$ in $\R^d$ such that
$$
\E_x\left[e^{i\xi\cdot(Y_t-Y_0)}\right]=e^{-t\psi(\xi)},
\quad \quad \mbox{ for every } x\in \R^d \mbox{ and } \xi\in \R^d,
$$
with
\begin{equation}\label{e:psi}
\psi(\xi)=\int_{\R^d}(1-e^{i\xi \cdot x}+i\xi \cdot x1_{\{|x|<1\}})
\nu^Y(x)dx, \qquad  \xi\in \R^d,
\end{equation}
where the L\'evy density $\nu^Y$ for $Y$ is given as
$$
\nu^Y(x):={\cal A} (d, - \alpha)|x|^{-d+\alpha}1_{\{|x|<1\}}, \qquad
x\in \R^d.
$$

Since $\psi(\xi)$ behaves like $|\xi|^{\alpha}$ near infinity
(see page 142 of \cite{KS}), by Proposition 28.1 in \cite{Sa} we know that the
process $Y$ has a smooth density $p^Y(t, x, y)$. Since $\psi(\xi)$
behaves like $|\xi|^2$ near the origin (again see page 142 of \cite{KS}),
it follows from Corollary 37.6 of \cite{Sa} that $Y$
is recurrent when $d=2$ and transient when $d\ge 3$.

For any open set $D$, we use $\tau^X_D$ to denote the first exit
time of $D$ for $X$, i.e., $\tau^X_D=\inf\{t>0: \, X_t\notin D\}$.
Given  an open set $D\subset \R^d$, we define
$X^D_t(\omega)=X_t(\omega)$ if $t< \tau^X_D(\omega)$ and
$X^D_t(\omega)=\partial$ if $t\geq  \tau^X_D(\omega)$, where
$\partial$ is a cemetery state. The process $X^D$ is usually called
the killed symmetric $\alpha$-stable process in $D$.

Similarly, we use $\tau^Y_D$
to denote the first exit time of $D$ for $Y$
and $Y^D$ to denote the process obtained by killing the process $Y$
upon exiting $D$.

Before we state more properties of truncated symmetric
$\alpha$-stable processes, we recall the following definitions.

\begin{defn}\label{def:har1} Let $D$ be an open subset of $\R^d$.
A locally integrable function $u$ defined on $\R^d$ taking values in
$(-\infty, \, \infty)$ and satisfying the condition $ \int_{\{x \in
\R^d ; |x| >1\}}|u(x)| |x|^{-(d+\alpha)} dx  <\infty $ is said to be

\begin{description}
\item{(1)}  harmonic for $X$ in $D$ if $$
\E_x\left[|u(X_{\tau^X_{B}})|\right] <\infty \quad \hbox{ and } \quad
u(x)= \E_x\left[u(X_{\tau^X_{B}})\right], \qquad x\in B, $$ for every
open set $B$ whose closure is a compact subset of $D$;
\item{(2}) regular  harmonic for $X$ in $D$ if
it is harmonic for $X$ in $D$  and
for each $x \in D$, $$ u(x)= \E_x\left[u(X_{\tau^X_{D}})\right]; $$
\item{(3}) harmonic for $X^D$ if it is harmonic for $X$ in $D$ and
vanishes outside $D$.
\end{description}
The corresponding concepts for
$Y$ can be defined similarly.
\end{defn}

In \cite{KS}, we have proved the following Harnack inequality for $Y$.

\begin{thm}[Theorem 4.9 in \cite{KS}]\label{T:Har0}
There exists $r_0 \in (0, \frac14)$ such that if  $r < r_0$ and
$x_{1}, x_{2}\in \R^d$ satisfy
 $|x_{1}-x_{2}|< Lr$
for some $L \le \frac{1}{r} - \frac12$,
then there exists a constant $c >0 $ depending only
on $d$ and $\alpha$, such that
$$
c^{-1}L^{-(d+\alpha)}u(x_{2})\,\leq \,u(x_{1})
\,\leq\, cL^{d+\alpha}u(x_{2})\,
$$
for every nonnegative function $u$
which is regular harmonic with respect to $Y$ in
$B(x_{1}, r)\cup B(x_{2},r)$.
\end{thm}

The Harnack inequality above  is similar to the Harnack inequality
for symmetric stable processes (Lemma 2 in \cite{B}), the difference
is that now we have to require that the two balls are not too far
apart. Because truncated stable processes can only make jumps of
size less than 1, one can easily see that, without the assumption
above, the Harnack inequality fails.

In this paper, we will mostly use the following simpler form of the
Harnack inequality, which is a direct consequence of Theorem
\ref{T:Har0}. From now on, $r_0$ will stand for the constant in
Theorem \ref{T:Har0}.

\begin{thm}\label{T:Har}
Suppose that $r\le r_0$. There exists a constant $c=
c(d, \alpha) >0$ such that
$$
c^{-1}\,u(y)\,\le\, u(x)\,\le\, c\,u(y), \qquad y\in B(x, \frac{r}2)
$$
for any nonnegative function $u$ which is
harmonic in $B(x, r)$ with respect to $Y$.
\end{thm}

\section{Green Function Estimates in Bounded Roughly
Connected $\kappa$-fat Open Sets}

Recall that $p^Y(t, x, y)=p^Y(t, x-y)$ is the transition density of
$Y$. For any bounded open set $D$ in $\R^d$,  let
$$
k_D(t,x,y)\,:=\,\E_x\left[p^Y(t-\tau^Y_D, Y_{\tau^Y_D},y):
\, \tau^Y_D < t \right]
$$
and
$$
p^Y_D(t,x,y):= p^Y(t,x,y)-k_D(t,x,y).
$$
It is well known that $p_D^Y(t,x,y)$ is the transition density of
$Y^D$.

Let $p^X(t,x,y)=p^X(t,x-y)$ be the transition density
function of $X$. It is well-known that
$$
p^X(t,x,y) \,\le\, c\, \left(t^{-d/\alpha} \wedge
\frac{t}{|x-y|^{d+\alpha}}\right), \quad (t,x,y) \in (0, \infty)
\times \R^d \times \R^d,
$$
for some $c=c(d, \alpha)>0$. Using the inequality above and the
smoothness of $p^Y(t,x,y)$, it is easy to see from Lemma 2.6 in
\cite{GR} that $p^Y(t,x)$ is bounded on the set $\{(t,x) : t>0, |x|
> \eps\}$ for $\eps >0$.

Let $C_0(\R^d)$ be the class of bounded continuous functions $f$ on
$\R^d$ with $\lim_{|x|\to\infty}f(x)=0$. We say a Markov process $Z$
in $\R^d$ has the Feller property if for every $g \in C_0(\R^d)$,
$\E_x[g(Z_t)]$ is in $C_0(\R^d)$ and $\lim_{t\to 0}\E_x[g(Z_t)]=g(x)$.
Any L\'evy process in  $\R^d$ has
the Feller property (for example, see \cite{Sa}). Using the
boundedness, the smoothness of $p^Y(t,x,y)$ and the separation
property for Feller processes (see Exercise 2 on page 73 of
\cite{CW}), it is routine (see, for instance, the proof of Theorem
2.4 \cite{CZ} ) to show that $Y^D$ has a jointly continuous and
symmetric density $p^Y_D(t, x, y)$. Moreover, from this one can
easily show that the Green function $G^Y_D(x,y):=\int_0^\infty
p^Y_D(t, x, y)dt$ of $Y^D$ is symmetric and continuous on $(D\times
D)\setminus\{(x, x): x\in D\}$.
Furthermore, for any $x\in D$, $G^Y_D(x, \cdot)=G^Y_D(\cdot, x)$
is harmonic for $Y$ in $D\setminus\{x\}$ and regular harmonic for $Y$ in $D\setminus B(x, \eps)$, $\eps>0$.

Let $G^X$ be the Green function of $X$ and
$G^X_{D}$ the Green function of
$X^{D}$. By Corollary 3.2 in \cite{GR} and the explicit
formula for $G^X$, we have that
for every bounded open set $D$ there exist
 $c_1=c_1(D, \alpha)>0 $ and $c_2=c_2(D, \alpha)>0 $ such that
\begin{equation}\label{e:G_ub}
G^Y_{D}(x, y)\,\le\, c_1\, G^X_{D}(x, y)\,\le\, c_2\, |x-y|^{-d+\alpha},
\quad \mbox{ for all } x, y \in D.
\end{equation}

Unlike the symmetric stable process $X$, the process $Y$ can only
make jumps of size less than 1. In order to guarantee $p^Y_D$ to be
strictly positive, we need to put the following assumption on $D$.

\begin{defn}\label{d:rc}
We say that an open set $D$ in $\R^d$ is roughly connected if for
every $x, y \in D$, there exist finite distinct connected components
$U_{1}, \cdots, U_{m}$ of  $D$ such that $ x \in U_{1}$, $y \in
U_{m}$ and dist$(U_{k}, U_{{k+1}}) <1$ for $1 \le k \le m-1$.
\end{defn}

The following result is proved in \cite{KS6}.

\begin{prop}[Proposition 4.4 in \cite{KS6}]\label{p:spdp}
For any bounded roughly connected open set $D$ in $\R^d$, the
transition density $p^Y_D(t,x,y)$ of $Y^D$ is strictly positive in
$(0, \infty)\times D \times D$.
\end{prop}

The process $X$ has a L\'evy system $(N, H)$ with $N(x, dy)={\cal
A}(d, -\alpha)|x-y|^{-(d+\alpha)}dy$ and $H_t=t$ (see \cite{FOT}).
Thus  for any open subset $D$ and $B \subset \R^d \setminus
\overline{D}$,
\begin{equation}\label{s_Levy}
\P_x\left(X_{\tau^X_D}  \in B \right) =\int_{B } K^X_D(x, z)
dz,\quad x \in D,
\end{equation}
where
\begin{equation}\label{e:pkX}
K^X_D(x, z):=
{\cal A}(d, -\alpha)\int_{D}
\frac{G^X_D(x,y)}{|y-z|^{d+\alpha}} dy.
\end{equation}
And the process $Y$ has a L\'evy system $(N^Y, H^Y)$ with $N^Y(x,
dy)={\cal A}(d, -\alpha)|x-y|^{-(d+\alpha)}1_{|x-y|<1}dy$ and
$H^Y_t=t$ (see \cite{FOT}).  Thus  for any open subset $D$ and $B
\subset \R^d \setminus \overline{D}$,
\begin{equation}\label{t_Levy}
\P_x\left(Y_{\tau^Y_D}  \in B \right) =\int_{B } K^Y_D(x, z)
dz,\quad x \in D,
\end{equation}
where
\begin{equation}\label{e:pkY}
K^Y_D(x, z):=
{\cal A}(d, -\alpha)\int_{D}
\frac{G^Y_D(x,y)}{|y-z|^{d+\alpha}} 1_{\{ |y-z| <1\}}dy.
\end{equation}
By \eqref{e:G_ub}, \eqref{e:pkY} and the continuity of Green function of $Y^D$,
$K^Y_D(\cdot, z)$ is continuous on $D$. Thus, using Lemma 4.2 in \cite{KS} and the explicit  formula for  $K^X_{B(w,r)}(x, z)$ (see \eqref{P_f} below), it is easy to see that the positive harmonic functions of $Y$ in any arbitrary open set $D$ are continuous in $D$.

Before we proceed, we recall the definition of
$\kappa$-fat set from \cite{SW}.

\begin{defn}\label{fat}
Let $\kappa \in (0,1/2]$. We say that an open set $D$ in
$\R^d$ is $\kappa$-fat if there exists $R>0$ such that for each
$Q \in \partial D$ and $r \in (0, R)$,
$D \cap B(Q,r)$ contains a ball $B(A_r(Q),\kappa r)$.
The pair $(R, \kappa)$ is called the characteristics of
the $\kappa$-fat open set $D$.
\end{defn}

Note that all Lipschitz domains and all non-tangentially accessible
domains (see \cite{JK} for the definition) are $\kappa$-fat.
Moreover, every  John domain is  $\kappa$-fat (see Lemma 6.3 in
\cite{MV}). Bounded $\kappa$-fat open sets may be disconnected.

Recently in \cite{KS4, KS5}, we have extended the concept of
intrinsic ultracontractivity to non-symmetric semigroups and proved
that for a large class of non-symmetric diffusions with
measure-valued drifts and potentials, the intrinsic
ultracontractivity is true in bounded domains under very mild
assumption. Using ideas in \cite{Ku}, we also showed in \cite{KS6}
that for a large class of non-symmetric L\'{e}vy processes with
discontinuous sample paths, the intrinsic ultracontractivity is true
in some bounded open sets. As a particular case of the general
result in \cite{KS6} we get that if $D$ is a roughly connected
bounded $\kappa$-fat open set, the semigroup of $Y^D$ is intrinsic
ultracontractive and the main results of \cite{KS6} are true (see
Example 4.2 in \cite{KS6}). In particular, the following is true.

\begin{prop}[Corollary 3.11 in \cite{KS6}]\label{p:lbG}
If $D$ is a bounded roughly connected $\kappa$-fat open set in
$\R^d$, then there exists constant $c=c(D, \alpha)>0$ such that
\begin{equation}\label{lbG}
c\, \E_x [ \tau^Y_D] \, \E_y [\tau^Y_D] \le\,
G^Y_D(x,y), \qquad (x,y) \in D \times D.
\end{equation}
\end{prop}

According to \cite{GR}, in the case when $D$ is a bounded Lipschitz
domain, the result above is also proved in \cite{G}.

In the remainder of this paper we assume that $D$ is a bounded
roughly connected $\kappa$-fat set with the characteristics  $(R,
\kappa)$. Without loss of generality, we assume $R <1$.

Before we further discuss properties of the Green function $G^Y_D$,
we recall some notations from \cite{KL}. Let $M:=2\kappa^{-1} $ and
fix $z_0 \in D$ with $2R/M < \rho_D(z_0) < R$ and let $\eps_1:= R
/(12M)$. For $x,y \in D$, we define $r(x,y): = \rho_D(x) \vee
\rho_D(y)\vee |x-y|$ and
$$
\BB(x,y):=\left\{ A \in D:\, \rho_D(A) > \frac1{M}r(x,y), \,  |x-A|\vee
|y-A| < 5 r(x,y)  \right\}
$$
if $r(x,y) <\eps_1 $, and $\BB(x,y):=\{z_0 \}$ otherwise.

It is well known (see \cite{CS1, Ku2}) that there exists a positive
constant $C_0$ such that $G^X_D(x,y) \le C_0 |x-y|^{-d+\alpha}$ for
$x,y \in D$ and that $ G^X_D(x,y) \ge G^X_{B(y, \rho _D(y))}(x,y)
\ge C_0^{-1} |x-y|^{-d+\alpha}$ if $|x-y| \le \rho_D(y)/2 $. Let
$C_1:=C_0 2^{d-\alpha} \rho_D(z_0)^{-d+\alpha}$ so that
$G^X_D(\cdot, z_0)$ is bounded from above by $C_1$ on $ D \setminus
B(z_0, \rho_D(z_0)/2)$. Now we define
$$
g(x):=  G^X_D(x, z_0) \wedge C_1.
$$
Note that if $\rho_D(z) \le 6 \eps_1$, then $|z-z_0| \ge \rho_D(z_0)
- 6 \eps_1 \ge \rho_D(z_0) /2$ since $6\eps_1<\rho_D(z_0)/4$, and
therefore $g(z )= G^X_D(z, z_0)$.

Using the Harnack inequality (Lemma 2 in \cite{B}) and the boundary
Harnack principle (Theorem 3.1 in \cite{SW}), the following Green
function estimates have been established by several authors. (See
Theorem 2.4 in \cite{H} and Theorem 1 in \cite{J}. Also see
\cite{KS2} for the case of non-symmetric diffusions.)

\begin{thm}[Theorem 2.4 in \cite{H}]\label{t:Gest}
Suppose that $D$ is a bounded $\kappa$-fat open set in $\R^d$. There
exists $c=c(D, \alpha)>0$ such that for every $x, y \in D$
\begin{equation}\label{e:Gest}
c^{-1}\,\frac{g(x) g(y)}{g(A)^2} \,|x-y|^{-d+\alpha} \,\le\,
G^X_D(x,y) \,\le\, c\,\frac{g(x) g(y)}{g(A)^2}\, |x-y|^{-d+\alpha},
\quad  A \in \BB(x,y).
\end{equation}
\end{thm}

\begin{lemma}\label{G:g0} Suppose that
$D$ is a bounded $\kappa$-fat open set in $\R^d$. There exists
$c=c(D, \alpha)>0$ such that
$$
c^{-1} \,\E_{x}\left[\tau^X_D\right] \,\le \,g(x)\, \le\,
c\,\E_{x}\left[\tau^X_D\right] , \qquad x \in D.
$$
\end{lemma}

\pf Choose a bounded open set $A \subset \R^d \setminus
\overline{D}$ with dist$(A,D) \ge 5 \eps_1$ and let
$h(x):=\P_x(X_{\tau^X_D}  \in A)$. By (\ref{s_Levy})-(\ref{e:pkX}),
there exists $c_1=c_1(\eps_1, A, D, \alpha)> 0$ such that
\begin{equation}\label{e:g01}
c_1^{-1} \E_{x}\left[\tau^X_D\right]   \,\le\, h(x)  \,\le\, c_1
\E_{x}\left[\tau^X_D\right] , \qquad x\in D.
\end{equation}
By the Harnack inequality for $X$ there exists
$c_2=c_2(\eps_1, A, D, \alpha)> 0$ such that
\begin{equation}\label{e:g02}
c_2^{-1} g(x)   \,\le\, h(x)  \,\le\, c_2 g(x) ,
\quad x \in \{y \in D; \rho_D (y)\ge \kappa \eps_1 \}.
\end{equation}
Recall that if $\rho_D(x) \le 6 \eps_1$, then $g(x)= G^X_D(x, z_0)$.
By (\ref{e:g01})-(\ref{e:g02}) it is enough to show that
\begin{equation}\label{e:g03}
c_3^{-1} G^X_D(x, z_0)  \,\le\, h(x)  \,\le\,
 c_3 G^X_D(x, z_0),  \quad x \in \{y \in D; \rho_D (y) < \kappa \eps_1 \}
\end{equation}
for some constant $c_3>0$. Since $2R/M < \rho_D(z_0) < R$ and
$\eps_1= R /(12M)$, for each $Q \in \partial D$, $G^X_D(x, z_0)$ and
$h(x)$ are regular harmonic in $D \cap B(Q,4 \eps_1)$, vanishing on
$D^c \cap B(Q,4 \eps_1)$. Thus by the boundary Harnack principle for
$X$ (Theorem 3.1 in \cite{SW}), we have for $x \in
D\cap B(Q,\eps_1 )$,
$$ c_4^{-1}
\kappa^{d+ \alpha}
\frac{G^X_D(A_{\eps_1}(Q),z_0)}{h(A_{\eps_1}(Q))}\, \le\,
\frac{G^X_D(x,z_0)}{h(x)} \,\le\, c_4 \,\kappa^{-d- \alpha}\frac{
G^X_D(A_{\eps_1}(Q),z_0)}{h(A_{\eps_1}(Q))}.
$$
Now applying (\ref{e:g02}) (note that
$\kappa \eps_1 \le \rho_D(A_{\eps_1}(Q)) < \eps_1 $) to the above,
we arrive at the conclusion of our lemma.\qed

\begin{lemma}\label{G:g2} Suppose
that $D$ is a bounded roughly connected $\kappa$-fat open set in
$\R^d$. Let $\eps >0$, then there exists $c=c(\eps, D, \alpha)>0$ such that for every $x,
y \in D$ satisfying $|x-y| \ge \eps$,
$$
G^X_D(x,y) \,\le\, c\,G^Y_D(x,y).
$$
\end{lemma}

\pf It follows from the Harnack inequality for $X$ (Lemma 2 in
\cite{B}) that there exists $c_1=c_1(D, \alpha)>0$ such that for
every $x, y \in D$,
$$
c_1^{-1} \,g(A_1) \,\le \,g(A_2)\, \le\, c_1\, g(A_1) , \qquad A_1,
A_2 \in \BB(x,y).
$$
Combining this with Theorem \ref{t:Gest} and Lemma \ref{G:g0}, we
get that there exists $c_2=c_2(D, \alpha)>0$ such that for every $x,
y \in D$,
$$
G^X_D(x,y) \,\le\, c_2\,\frac{\E_{x}[\tau^X_D]
\E_{y}[\tau^X_D]}{g^2(A)|x-y|^{d-\alpha}}, \qquad A \in \BB(x,y).
$$
Note that if $|x-y| \ge \eps$, there exists $c_3=c_3(D, \alpha, \eps)>0$ such that $g(A) >c_3$ for $A \in \BB(x,y)$ because either $\rho_D(A) > \eps/M$ or $A=z_0$.  Moreover,
it follows from Lemma 2.4 of \cite{GR} that there exist a constant
$c_4=c_4(D, \alpha)>0$ such that
$$
c_4^{-1}\,\E_{x}\left[\tau^X_D\right]\, \le\,
\E_{x}\left[\tau^Y_D\right]\, \le\, c_4\,
\E_{x}\left[\tau^X_D\right], \quad x\in D.
$$
Therefore, by Proposition
\ref{p:lbG} we get
$$
G^X_D(x,y) \,\le\, c_2\,c_3^{-2} \,c_4^2\,\eps^{\alpha-d}\,\E_{x}[\tau^Y_D]
\E_{y}[\tau^Y_D]   \,\le\, c_5\,G^Y_D(x,y)
$$
for some positive constant $c_5=c_5(D, \alpha, \eps)>0$
\qed

\begin{thm}\label{t:ge}
Suppose that $D$ is a bounded roughly connected $\kappa$-fat open
set in $\R^d$. Then there exists $c=c(D, \alpha) >0$ such that
\begin{equation}\label{e:ge}
c^{-1} \,G^Y_D(x,w) \,\le\, G^X_D(x,w) \,\le\, c\, G^Y_D(x,w), \quad
(x,w) \in D \times D.
\end{equation}
\end{thm}

\pf By (\ref{e:G_ub}) and Lemma \ref{G:g2}, we only need to show the
second inequality in (\ref{e:ge}) for $|x-w|< \eps$, where $\eps>0$
is to be chosen later. By Corollary 3.6 in \cite{GR} (also see (11)
in \cite{GR}),

\begin{eqnarray*}
G^X_D(x,w) &\le &  G^Y_D(x,w) \,+\,{\cal A}(d, -\alpha)
\int_D \int_{D}  G^X_D(x,y) G^X_D(z,w) 1_{\{|y-z|>1\}}(y,z)|y-z|^{-(d+\alpha)}
dydz\\
&\le &  G^Y_D(x,w) \,+\,{\cal A}(d, -\alpha)
\int_D \int_{D}  G^X_D(x,y) G^X_D(z,w)
dydz \\
&=&  G^Y_D(x,w) \,+\,{\cal A}(d, -\alpha) \,\E_x[\tau^X_{D}] \,\E_w[\tau^X_{D}].
\end{eqnarray*}
Applying Theorem \ref{t:Gest} and Lemma \ref{G:g0}, we get
$$
G^X_D(x,w) \, \le \,  G^Y_D(x,w) \,+\, c_1 \, G^X_D(x,w) |x-w|^{d-\alpha}
$$
for some positive constant $c_1>0$.
Choose $\eps>0$ small so that
$
2 c_1 G^X_D(x,w) |x-w|^{d-\alpha}  \le G^X_D(x,w)$ if  $|x-w| < \eps$.
Thus
$
 G^X_D(x,w) \le  2  G^Y_D(x,w)$,  if $ |x-w| < \eps.
$
\qed

Let $z_0 \in D$, $C_1$ and $\BB(x,y)$ be the same as defined before
Theorem \ref{t:Gest}, and  let
$$
g^Y(x):=  G^Y_D(x, z_0) \wedge C_1.
$$
By  Theorem \ref{t:ge}, $g^Y$ is comparable to $g$. Thus we can
easily get the following Green function estimates for Y in bounded
roughly connected $\kappa$-fat open sets from Theorem \ref{t:Gest}
and Theorem \ref{t:ge}.

\begin{thm}\label{t:GestY}
Suppose that $D$ is a bounded roughly connected $\kappa$-fat open
set in $\R^d$. Then there exists $c=c(D, \alpha)>0$ such that for
every $x, y \in D$
\begin{equation}\label{e:GestY}
c^{-1}\,\frac{g^Y(x) g^Y(y)}{g^Y(A)^2} \,|x-y|^{-d+\alpha} \,\le\,
G^Y_D(x,y) \,\le\, c\,\frac{g^Y(x) g^Y(y)}{g^Y(A)^2}\, |x-y|^{-d+\alpha},
\quad  A \in \BB(x,y).
\end{equation}
\end{thm}

Combining Theorem \ref{t:ge} and the generalized 3G theorem for $X$
(Theorem 1.1 in \cite{KL}),
we get the generalized 3G theorem for truncated stable
processes on roughly connected
$\kappa$-fat open sets.

\begin{thm}[Generalized 3G theorem for $Y$]\label{3GY}
Suppose $D$ is a bounded roughly connected $\kappa$-fat open set in
$\R^d$. Then there exist positive constants $c=c(D, \alpha)$ and $
\gamma <\alpha$ such that for every $x, y, z, w \in D$
\begin{equation} \label{e:3GY}
\frac{G^Y_D(x,y) G^Y_{D}(z,w)} { G^Y_{D}(x,w)} \le c
\left(\frac{|x-w|\wedge |y-z|}{|x-y|  } \vee 1 \right)^\gamma
\left(\frac{|x-w|\wedge |y-z|}{|z-w|} \vee 1 \right)^\gamma
\frac{|x-w|^{d-\alpha}} {|x-y|^{d-\alpha} |z-w|^{d-\alpha}}.
\end{equation}
\end{thm}

A  bounded open set $U$ in $\R^d$ is said to be a $C^{1,1}$ open set
if there is a localization radius $R_0>0$ and a constant $\Lambda
>0$ such that for every $Q\in \partial U$, there is a
$C^{1,1}$-function $\phi=\phi_Q: \R^{d-1}\to \R$ satisfying $\phi
(0) = |\nabla\phi (0)|=0$, $\| \nabla \phi  \|_\infty \leq \Lambda$,
$| \nabla \phi (x)-\nabla \phi (z)| \leq \Lambda |x-z|$, and an
orthonormal coordinate system $y=(y_1, \cdots, y_{d-1},
y_d):=(\tilde y, y_d)$ such that $ B(Q, R_0)\cap U=B(Q, R_0)\cap \{
y: y_d > \phi (\tilde y) \}$.

Recently, sharp estimates (even in terms of $\alpha$ and $d$) on
$G^X_D$ for a bounded $C^{1,1}$ open set $D$ were obtained in
\cite{C1, CS4} (see also \cite{CS1, Ku2} for estimates on $G^X_D$
for bounded $C^{1,1}$ domains). Thus by Theorem \ref{t:ge}, we have
the sharp estimates on $G^Y_D$.

 \begin{thm}\label{t:geC}
Suppose $D$ is a bounded roughly connected $C^{1,1}$ open set in
$\R^d$. Then there exists $c=c(D, \alpha)>0$ such that for every $x,
y \in D$
\begin{equation}\label{eqn:G1}
c^{-1} \left( \frac1{|x-y|^{d-\alpha}} \wedge  \frac{\rho_D
(x)^{\alpha/2} \rho_D (y)^{\alpha/2}} {|x-y|^d}\right) \,\le\,
G^Y_D(x, y) \,\le\,
 c\left( \frac1{|x-y|^{d-\alpha}} \wedge  \frac{\rho_D (x)^{\alpha/2}
\rho_D (y)^{\alpha/2}} {|x-y|^d}\right).
\end{equation}
\end{thm}

\section{Martin Boundary and Martin Representation}\label{sec-mbr}

In this section we will always assume that $D$ is a bounded roughly
connected $\kappa$-fat open set in $\R^d$ with the characteristics
$(R, \kappa)$. We are going to apply Theorem \ref{t:ge} to study the
Martin boundary of $D$ with respect to $Y$. The argument in this
section is similar to that of Section 4 of \cite{SW}.

We recall from Definition \ref{fat} that  for each $Q \in \partial
D$ and $r \in (0, R)$, $A_r(Q)$ is a point in $D \cap B(Q,r)$
satisfying $B(A_r(Q),\kappa r)  \subset D \cap B(Q,r)$. Combining
the boundary Harnack principle for $X$ (Theorem 3.1 in \cite{SW})
and Theorem \ref{t:ge}, we get the following boundary Harnack
principle for  Green functions of $Y$ which will play an important
role in this section.

\begin{thm}\label{BHP2}
There exists a constant $c=c(D, \alpha)
>1$ such that for any $Q \in \partial D$, $r \in (0,R)$
and $z,w \in D \setminus B(Q,2r)$, we have $$ c^{-1}\,
\frac{G_D^Y(z, A_r(Q))}{G_D^Y(w, A_r(Q))} \,\le\,
\frac{G_D^Y(z,x)}{G_D^Y(w, x)} \,\le\, c\, \frac{G_D^Y(z,A_r(Q))}
{G_D^Y(w, A_r(Q))} ,\quad x\in D\cap B\left(Q,\frac{r}{2}\right). $$
\end{thm}

Recall that $r_0>0$ is the constant from Theorem \ref{T:Har0}.  Without
loss of generality, we will assume $r_0 < R$.
Recall that $A(x, a,b):=\{ y \in \R^d: a \le |y-x| <b  \}.$

The following  result was proved in \cite{KS} for harmonic functions
of $Y$ when $D$ is bounded convex domain, which is analogous to
Lemma 5 of \cite{B}. It was observed in \cite{KL} (Lemma 3.1 in
\cite{KL}) that it is valid for a large class of jump processes. We
reproduce the proof here for the sake of completeness.

\begin{lemma}\label{l:5B}
There exist positive constants $c=c(D, \alpha)$ and
$\gamma=\gamma(D, \alpha)< \alpha$ such that for any $Q\in \partial
D$ and $r\in (0, r_0)$, and nonnegative function $u$ which is
harmonic with respect to $Y$ in $D \cap B(Q, r)$ we have
\begin{equation}\label{e:gamma}
u(A_s(Q))\,\ge\, c\,(s/r)^{\gamma}\,u(A_r(Q)), \qquad s\in (0, r).
\end{equation}
\end{lemma}

\pf
Without loss of generality, we may assume $Q=0$.
Fix $r <  r_0$ and
let
$$
 \eta_k\,:=\,\left(\frac{\kappa}2\right)^{k}r,
 \quad A_k\,:=\, A_{ \eta_k}(0)
\quad \mbox{ and } \quad B_k\,:=\,B(A_k,  \eta_{k+1}), \quad k=0,1, \cdots.
$$
Note that the $B_k$'s are disjoint. So by the harmonicity of $u$, we have
$$
u(A_k)
\,\ge\, \sum_{l=0}^{k-1} \E_{A_k}\left[u(Y_{\tau^Y_{B_k}}):\,
Y_{\tau^Y_{B_k}} \in B_l \right]\\
\,=\, \sum_{l=0}^{k-1} \int_{B_l} K_{B_k}^Y(A_k, z) u(z) dz.
$$
Since $r < r_0 < \frac14$, we have by Lemma 4.2 of \cite{KS} that
for every $w \in \R^d$ and $z \in A(w, r, 1-r)$,
\begin{equation}\label{CC1}
K^X_{B(w,r)}(x,z) \,\le\,  K_{B(w,r)}^Y(x,z),\qquad x \in
B(w,r).
\end{equation}
 Now (\ref{CC1}) and
Theorem  \ref{T:Har} imply that
$$
 \int_{B_l} K_{B_k}^Y(A_k, z) u(z) dz \,\ge\, c_1\, u(A_l)
\int_{B_l} K^X_{B_k}(A_k, z) dz
$$
for some constant $c_1=c_1(d, \alpha)>0$. It is well-known that
\begin{equation}\label{P_f}
K^X_{B(w,r)}(x,z)\,=\,c_2\,
\frac{(r^2-|x-w|^2)^{\frac{\alpha}2}}{(|z-w|^2-r^2)^{\frac{\alpha}2}}
\frac1{|x-z|^d}
\end{equation}
for some constant $c_2=c_2(d, \alpha) > 0$. Using  (\ref{P_f}), one
can easily check that
$$
\int_{B_l}K^X_{B_k}(A_k, z)dz \,\ge\, c_3 \,
\left(\frac{\kappa}{2}\right)^{(k-l) \alpha}, \quad z \in B_l,
$$
for some constant $c_3=c_3(d, \alpha)>0$. Therefore,
$$
\left( \eta_k\right)^{-\alpha} u(A_k) \,\ge \, c_4 \sum_{l=0}^{k-1}
\left( \eta_l\right)^{-\alpha} u(A_l)
$$
for some constant $c_4=c_4(d, \alpha)>0$. Let $a_k :=
 (\eta_k)^{-\alpha}u(A_k)$ so that $a_k \ge  c_4\sum_{l=0}^{k-1}  a_l$. By
induction, one can easily check that $ a_k  \ge c_5 (1+c_4/2)^{k} $
for some constant $c_5=c_5(d, \alpha)>0$. Thus, with $ \gamma =
\alpha - {\ln(1+\frac{c_4}2)} (\ln (2/\kappa))^{-1}, $
(\ref{e:gamma}) is true for $s=  \eta_k$. For the other values, we apply
the Harnack inequality (Theorem \ref{T:Har}). \qed

The next lemma is the Green function version of
Lemma 5.8 of \cite{KS}. The novelty here is that
we are now dealing with bounded roughly connected $\kappa$-fat open sets
rather than bounded convex domains.

\begin{lemma}\label{l:la}
Suppose $Q \in \partial D$ and $r \in (0, R)$. If $w \in D\setminus
B(Q, r)$, then
\[
G_D^Y(A_r(Q), w)\, \ge\, c\, r^{\alpha} \int_{A(Q, r, 1+\frac{r}2)}
|z-Q|^{-d-\alpha}G_D^Y(z, w)dz
\]
for some constant $c=c(D, \alpha)>0$.
\end{lemma}

\pf It follows from Theorem \ref{t:ge} that it is enough to prove
this lemma for $G^X_D$ instead of $G^Y_D$. Without loss of
generality, we may assume $Q=0$. Fix $w \in D\setminus B(0, r)$ and
let $A:=A_r(0)$ and $u(\cdot) := G_D^X(\cdot, w)$. Since $u$ is
regular harmonic in $D\cap B(0, (1-\kappa/2)r)$ with respect to $X$,
by (\ref{e:pkX}) we have
\begin{eqnarray*}
&&u(A) \,\ge\, \E_A \left[ u\left( X_{\tau^X_{D \cap
B(0,(1-\kappa/2)r)}}\right); X_{\tau^X_{D \cap B(0,(1-\kappa/2)r)}}
\in
A(0, r, 1+\frac{r}2) \right]\\
&&= \int_{A(0,r , 1+\frac{r}2)}K^X_{D\cap
B(0, (1-\kappa/2)r)}(A, z)u(z)dz\\
&&=\int_{A(0, r, 1+\frac{r}2)         }
{\cal A}(d, -\alpha)\int_{D\cap B(0,  (1-\kappa/2)r ) }
 \frac{G_{D\cap B(0,(1-\kappa/2)r )}^X(A,y)}{|y-z|^{d+\alpha}}
 dyu(z)dz.
\end{eqnarray*}
Since $B(A, \kappa r/2)\subset D \cap B(0, (1-\kappa/2)r )$,
by the monotonicity of
the Green functions,
$$
G_{D\cap B(0,(1-\kappa/2)r )}^X(A,y) \, \ge \,
G^X_{B(A, \kappa r/2)}(A,y),
\quad y \in B(A, \kappa r/2).
$$
Thus
$$
u(A) \ge  \int_{A(0,r , 1+\frac{r}2)   }
{\cal A}(d, -\alpha)\int_{ B(A, \kappa r/2)  }
\frac{G^X_{B(A, \kappa r/2)}(A,y)}{|y-z|^{d+\alpha}} dyu(z)dz
= \int_{ A(0, r, 1+\frac{r}2)   }
K^X_{B(A, \kappa r/2)}(A, z)u(z)dz,
$$
which is equal to
$$
c_1\int_{A(0, r, 1+\frac{r}2)  }
\frac{(\kappa r/2)^{\alpha}}{(|z-A|^2-
(\kappa r/2)^2)^{\frac{\alpha}2}}\frac1{|z-A|^d}u(z)dz
$$
for some positive constant $c_1=c_1(d, \alpha)$ by (\ref{P_f}). Note
that $|z-A|\le  2|z|$ for $z\in A(0,r , 1+\frac{r}2) $. Hence
\begin{equation}\label{e:22}
u(A)\,\ge\, c_2\,\kappa^{\alpha}\,r^{\alpha}\int_{A(0,r ,1+\frac{r}2 ) }
\frac{u(z)}{|z|^
{d+\alpha}}dz
\end{equation}
for some constant $c_2=c_2(d, \alpha)>0$. \qed

By taking $r_0$ smaller if necessary, we get from Lemma 4.2 of
\cite{KS} that for every $r\le r_0$, $w \in \R^d$ and
$z \in A(w, r,\infty)$, we have
\begin{equation}\label{CC2}
K_{B(w,r)}^Y(x,z) \,\le\, 2 K_{B(w,r)}^X(x,z),\qquad x \in B(w,r).
\end{equation}

The next lemma is the Green function version of Lemma 5.14 of \cite{KS}.

\begin{lemma}\label{l:14B}
There exist positive constants $c_1=c_1(D,\alpha)$  and
$c_2=c_2(D,\alpha)<1$ such that for any $Q \in \partial D $, $r\in
(0, r_0)$ and  $w \in D \setminus B(Q,4r)$,  we have
$$
\E_x\left[G_D^Y(Y_{\tau^Y_{D \cap B_k}}, w):\,Y_{\tau^Y_{D
\cap
B_k}} \in A(Q, r, 1+4^{-k}r)
 \right] \,\le\, c_1\,c_2^{k} \, G_D^Y(x,w), \quad x \in D \cap B_k,
$$
where $B_k:=B(Q, 4^{-k}r)$, $ k=0,1, \cdots$.
\end{lemma}

\pf Without loss of generality, we may assume $Q=0$. Fix $r <r_0$
and  $w \in D \setminus B(0,4r)$.  Let $\eta_k:=4^{-k}r $,
$B_k:=B(0,\eta_k)$ and
$$
u_k(x)  \,:=\, \E_x\left[G^Y_D(Y_{\tau^Y_{D \cap B_k}}, w);
Y_{\tau^Y_{D \cap B_k}}
 \in A(0, r, 1+\eta_k) \right], \quad x \in D \cap B_k.
$$
Note that
\begin{eqnarray*}
u_{k+1}(x) &=&  \E_x\left[G^Y_D(Y_{\tau^Y_{D \cap B_{k+1}}}, w);\,
Y_{\tau^Y_{D \cap B_{k+1}}}
 \in A(0, r, 1+\eta_{k+1}) \right] \\
&=&  \E_x\left[G^Y_D(Y_{\tau^Y_{D \cap B_{k+1}}}, w);\, \tau^Y_{D
\cap B_{k+1}} = \tau^Y_{D \cap B_{k}}, ~ Y_{\tau^Y_{D \cap B_{k+1}}}
 \in A(0, r, 1+\eta_{k+1})\right] \\
&=&  \E_x\left[ G^Y_D(Y_{\tau^Y_{D \cap B_k}}, w)  ;\, \tau^Y_{D
\cap B_{k+1}} = \tau^Y_{D\cap B_{k}}, ~ Y_{\tau^Y_{D \cap B_{k}}}
 \in A(0, r, 1+\eta_{k+1}) \right] \\
&\le& \E_x\left[G^Y_D(Y_{\tau^Y_{D \cap B_k}}, w);\, Y_{\tau^Y_{D
\cap B_{k}}}
 \in A(0, r, 1+\eta_{k+1}) \right]\\
 &\le &
 \E_x\left[G^Y_D(Y_{\tau^Y_{D \cap B_k}}, w);\, Y_{\tau^Y_{D \cap
B_{k}}}
 \in A(0, r, 1+\eta_{k}) \right].
\end{eqnarray*}
Thus
\begin{equation}\label{e:dec1}
u_{k+1}(x) \,\le\, u_{k}(x).
\end{equation}
Let
$A_k\,:=\,A_{\eta_k}(0) $.
We have
\begin{eqnarray*}
&&u_{k}(A_k) \,= \, \E_{A_k}\left[G^Y_D(Y_{\tau^Y_{D \cap B_k}}, w);\,
Y_{\tau^Y_{D \cap B_{k}}}
 \in A(0, r, 1+\eta_{k}) \right] \\
&&\le  \E_{A_k}\left[G^Y_D(Y_{\tau^Y_{ B_k}}, w)  ;\, Y_{\tau^Y_{
B_{k}}}
 \in A(0, r, 1+\eta_{k})\right]
\,\le\,  \int_{A(0, r, 1+\eta_{k})} K_{B_k}^Y(A_k,z) G^Y_D(z,w) dz.
\end{eqnarray*}
Since $r<r_0$, by  (\ref{CC2}) and (\ref{P_f}), we get that for $z
\in A(0, r, 1+\eta_{k})$,
$$
K_{B_k}^Y(A_k,z) \,\le\,2\,K_{B_k}^X (A_k,z) \,\le\,  c_1
\,\frac{4^{-k\alpha} r^\alpha}{|z|^{d+\alpha}}
$$
for some constant $c_1=c_1(d, \alpha)>0$ and $k=1,2, \cdots$.
Therefore
\begin{equation}\label{e:dec2}
 u_{k}(A_k) \,\le\, c_2\,4^{-k\alpha}  r^\alpha\int_{A(0, r, 1+\eta_{k})}
G^Y_D(z,w)\frac{dz}{|z|^{d+\alpha}}, \quad k=1,2, \cdots
\end{equation}
for some constant $c_2=c_2(d, \alpha)>0$. From Lemma \ref{l:la}, we
have
\begin{equation}\label{e:dec3}
 G^Y_D(A_0,w) \,\ge\, c_3\ r^\alpha\int_{A(0, r, 1+\frac{r}2)}
G^Y_D(z,w)\frac{dz}{|z|^{d+\alpha}}
\end{equation}
for some constant $c_3=c_3(D, \alpha)>0$. (\ref{e:dec2}) and
(\ref{e:dec3}) imply that $ u_{k}(A_k) \,\le\, c_4\,4^{-k\alpha}
G^Y_D(A_0,w)$ for some constant $c_4=c_4(D, \alpha)>0$. On the other
hand, using Lemma \ref{l:5B}, we get $G_D^Y(A_0,w) \,\le\,
c_5\,4^{k\gamma } G^Y_D(A_k,w)$ for some constant $c_5=c_5(D,
\alpha)>0$. Thus, $ u_{k}(A_k) \,\le\,
c_4c_5\,4^{-k(\alpha-\gamma)}G^Y_D(A_k,w)$. By Theorem \ref{BHP2},
we have
$$
\frac{u_k(x)}{G_D^Y(x, w)}    \,\le\,
\frac{u_{k-1}(x)}{G_D^Y(x, w)}
\,\le\, c_6\, \frac{u_{k-1}(A_{k-1})}{G_D^Y(A_{k-1}, w)}  \,\le\,
c_4c_5c_6\,4^{-(k-1)(\alpha-\gamma)}
$$
for some constant $c_6=c_6(D, \alpha)>0$ and $k=1,2, \cdots.$ \qed

Let $x_0\in D$ be fixed and set
\[
M^Y_D(x, y):=\frac{G^Y_D(x, y)}{G^Y_D(x_0, y)}, \qquad x, y\in D,~
y\neq x_0.
\]
$M^Y_D$ is called the Martin kernel of $D$ with respect to $Y$.

Now the next theorem follows from Lemma \ref{l:5B}, Theorem
\ref{BHP2} and Lemma \ref{l:14B} (instead of Lemma 5,  Lemma 13 and
Lemma 14 in \cite{B} respectively) in very much the same way as in
the case of symmetric stable processes in Lemma 16 of \cite{B} (with
Green functions instead of harmonic functions). We omit the details.

\begin{thm}\label{t2.2}
There exist positive constants $r_1$, $M_1$, $c$ and $\gamma$
depending on $D$ and $\alpha$ such that for any $Q \in \partial D $,
$r < r_1$ and $z \in D \setminus B(Q, M_1 r)$, we have
\[
\left|M^Y_D(z, x)-M^Y_D(z, y)\right| \,\le\,
c\,\left(\frac{|x-y|}r\right)^{\gamma}, \qquad x, y\in  D \cap B(Q,
r).
\]
In particular, the limit $\lim_{D \ni y\to w} M^Y_D(x, y)  $ exists
for every $w\in \partial  D$.
\end{thm}

There is a compactification $D^M$ of $D$, unique up to a homeomorphism,
such that $M^Y_D(x, y)$ has a continuous
extension to $D\times (D^M\setminus\{x_0\})$ and $M^Y_D(\cdot, z_1)
=M^Y_D(\cdot, z_2)$ if and only if $z_1=z_2$. (See, for instance,
 \cite{KW}.)
The set $\partial^MD=D^M\setminus D$ is called the Martin boundary
of $D$. For $z\in \partial^MD$, set $M^Y_D(\cdot, z)$ to be zero in
$D^c$.

A positive harmonic function $u$ for  $Y^D$ is minimal if, whenever
$v$ is a positive harmonic function for $Y^D$ with $v\le u$ on $D$,
one must have $u=cv$ for some constant $c$. The set of points $z\in
\partial^MD$ such that $M^Y_D(\cdot, z)$ is minimal harmonic
for $Y^D$ is called the minimal Martin boundary of $D$.

For each fixed $z\in \partial D$ and $x\in D$,
let
\[
M^Y_D(x, z):=\lim_{D\ni y\to z}M^Y_D(x,y),
\]
which exists by Theorem \ref{t2.2}. For each $z\in
\partial D$, set $M^Y_D(x, z)$
to be zero for $x\in D^c$.

\begin{lemma}\label{l:MH1}
For every $z \in \partial D$ and $B \subset \overline{B} \subset D$,
$M^Y_D(Y_{\tau^Y_{B}} , z)$ is $\P_x$-integrable.
\end{lemma}
\pf
Take a sequence $\{z_m\}_{m \ge 1} \subset D\setminus
\overline{B}$ converging to $z$.
Since $M^Y_D (\cdot , z_m )$ is regular harmonic for $Y$ in $B$,
by Fatou's lemma and Theorem \ref{t2.2},
$$  \E_x \left[ M^Y_D\left(Y_{\tau^Y_{B}} , z\right)\right]\,=\,
\E_x \left[ \lim_{m\to \infty} M^Y_D\left(Y_{\tau^Y_{B}} ,
z_m\right) \right]\,\leq\, \liminf_{m\to \infty} M^Y_D(x, z_m)\,=\,
M^Y_D(x,z) \,<\,\infty .
$$
\qed

\begin{lemma}\label{l:MH2}
For every $z \in \partial D$ and  $x \in D$,
\begin{equation}\label{eqn:4.1}
M^Y_D(x,z) \,=\, \E_x \left[M^Y_D\left(Y^D_{\tau^Y_{B(x,r)}} ,
z\right)\right], \quad \mbox{ for
 every } 0<r<r_0 \wedge \frac12\rho_D(x).
\end{equation}
\end{lemma}
\pf Fix $z \in \partial D$, $x \in D$ and $r<r_0 \wedge
\frac12\rho_D(x)<R$. Let $z_m:=A_{r/m}(z)$ for $m \ge 2$ so that
$$ B(z_m , \, \kappa r/m ) \,\subset \, B(z,\,  r/m
) \cap D
 \,\subset  \,B(z, \,  2r/m  )\cap D  \,\subset  \,B(z, r)\cap D
 \,\subset \, D \setminus B(x,r)
$$
for all $m \ge2 $. Thus by the harmonicity of $M^Y_D(\cdot , z_m )$,
we have
$$M^Y_D(x,z_m) \,=\, \E_x \left[M^Y_D\left(Y_{\tau^Y_{B(x,r)}} , z_m\right)\right].$$

On the other hand, by Theorem \ref{BHP2}, there exist constants $m_0
\ge 2$ and $c_1>0$ such that for every $w \in D \setminus B (z,
2r/m)$ and  $y \in D \cap B(z, r/(2m))$,
$$
M^Y_D(w, z_m) \,=\, \frac{G^Y_D(w, z_m)}{G^Y_D(x_0 , z_m )} \,\leq\,
c_1\, \frac{G^Y_D(w, y)}{G^Y_D(x_0 , y )} \,=\, c_1\, M^Y_D(w,y),
\quad m\geq m_0.
$$
Letting $y \rightarrow z\in \partial D$ we get
\begin{equation}\label{eqn:4.2}
M^Y_D(w, z_m) \,\leq\, c_1\, M^Y_D (w, z), \quad m\geq m_0,
\end{equation}
for every $w \in D \setminus B(z, 2 r/m).$

To prove (\ref{eqn:4.1}), it  suffices to show that
$\{ M^Y_D(Y_{\tau^Y_{B(x,r)}}, z_m ) : m\geq m_0 \}$  is $\P_x $-uniformly
integrable.
Since $M^Y_D(Y_{\tau^Y_{B(x,r)}}, z)$ is $\P_x$-integrable
by Lemma \ref{l:MH1}, for
any $\eps>0$, there is an $N_0>1$ such that
\begin{equation}\label{eqn:c4.3}
\E_x \left[M^Y_D\left(Y_{\tau^Y_{B(x,r)}}, z\right) ;\,
M^Y_D\left(Y_{\tau^Y_{B(x,r)}}, z\right) > N_0/c_1 \right] \,<\,
\frac{\eps}{4c_1}.
\end{equation}
Note that by (\ref{eqn:4.2}) and (\ref{eqn:c4.3})
\begin{eqnarray*}
&& \E_x \left[M^Y_D\left(Y_{\tau^Y_{B(x,r)}}, z_m\right) ; \,
M^Y_D\left(Y_{\tau^Y_{B(x,r)}}, z_m\right) > N_0 ~\mbox{ and }~
Y_{\tau^Y_{B(x,r)}} \in D \setminus B(z, 2r/m)\right] \\
&\leq & c_1\, \E_x \left[M^Y_D\left(Y_{\tau^Y_{B(x,r)}}, z\right) ;
\, c_1 M^Y_D\left(Y_{\tau^Y_{B(x,r)}}, z\right) > N_0 \right]
 \,<\,   c_1\,\frac{ \eps}{4c_1} \,=\, \frac{\eps}{4} .
\end{eqnarray*}
Since $r < r_0$, by \eqref{P_f}, \eqref{CC2}, we have for $m\geq
m_0$,
\begin{eqnarray}
&&\E_x \left[M^Y_D\left(Y^D_{\tau^Y_{B(x,r)}}, z_m\right) ; \,
Y_{\tau^Y_{B(x,r)}} \in D
 \cap B(z, 2r /m)\right]\nonumber\\
&&\leq \, 2\, \E_x \left[ M^Y_D\left(X^D_{\tau^X_{B(x,r)}}, z_m
\right) 1_{D \cap B(z, 2r/m)}
     \left(X^D_{\tau^X_{B(x,r)}}\right)\right]\nonumber\\
&&\leq \, c_2 \, \int_{B(z, 2r/m)} M^Y_D(w, z_m) dw \,=\, c_2 \,
G^Y_D(x_0, z_m)^{-1} \int_{B(z, 2r/m)} G^Y_D(w, z_m)
dw\label{e:new1}
\end{eqnarray}
for some $c_2>0$.
Note that, by Lemma \ref{l:5B}, there exist $c_3
>0$, $c_4
>0$ and $\gamma < \alpha$ such that
\begin{equation}\label{e:new2}
G^Y_D(x_0, z_m) \,\ge\, c_3\, \left(
\frac{r/m}{r/m_0}\right)^{\gamma}G^Y_D(x_0, z_{m_0}) \,\ge\,
c_4\,\frac1{m^{\gamma}}.
\end{equation}
It follows from (\ref{e:new1})-(\ref{e:new2}) and Theorem \ref{t:ge}
that there exist $c_5>0$ and $c_6>0$ such that
$$
\E_x \left[M^Y_D\left(Y^D_{\tau^Y_{B(x,r)}}, z_m\right) ; \,
Y_{\tau^Y_{B(x,r)}} \in D
 \cap B(z, 2r /m)\right] \,\le\,  c_5 \,m^{\gamma}
 \int_{B(z, 2r/m)} \frac{dw}{|w-z_m|^{d-\alpha}} \,\leq \, c_6 \,
\frac{1}{m^{\alpha-\gamma}}.
$$
Therefore by taking $N$ large enough we have for $m\geq N$,
\begin{eqnarray*}
&& \E_x  \left[M^Y_D\left(Y_{\tau^Y_{B(x,r)}}, z_m\right);\,
M^Y_D\left(Y_{\tau^Y_{B(x,r)}}, z_m\right) > N\right]\\
&\leq&  \E_x  \left[M^Y_D\left(Y_{\tau^Y_{B(x,r)}}, z_m\right) ; \,
Y_{\tau^Y_{B(x,r)}} \in D \cap B(z, 2r/m)\right]\\
&&+ \E_x  \left[M^Y_D\left(Y_{\tau^Y_{B(x,r)}}, z_m\right) ; \,
M^Y_D\left(Y_{\tau^Y_{B(x,r)}}, z_m\right) > N \,\mbox{ and }\,
Y_{\tau^Y_{B(x,r)}} \in D \setminus B(z, 2r/m)\right] \\
&<&  c_6 \, \frac{1}{m^{\alpha-\gamma}} \,+\, \frac{\eps}{4}
\,\,<\,\,  \eps.
\end{eqnarray*}
As each $M^Y_D(Y_{\tau^Y_{B(x,r)}}, z_m)$ is $\P_x$-integrable,
we conclude that
$\{ M^Y_D(Y_{\tau^Y_{B(x,r)}}, z_m ) : m\geq m_0 \}$  is uniformly integrable
under $\P_x$. \qed

The two lemmas above imply that $M^Y_D(\cdot, z)$ is harmonic.

\begin{thm}\label{T:L4.3}
For every $z \in \partial D$,
the  function $x\mapsto M^Y_D(\cdot, z)$ is harmonic in $D$
with respect to $Y$.
\end{thm}

\pf Fix $z\in \partial D$ and let $h(x):=M_D^Y(x,z)$. For any open
set $D_1 \subset\overline{D_1}\subset D$ we can always take a smooth
open set $D_2$ such that $D_1\subset\overline{D_1}\subset D_2
\subset\overline{D_2}\subset D  $. Thus by the strong Markov property,
it is enough to show that for any $x\in D_2$,
\[
h(x)= E_x\left[h\left(Y_{\tau^Y_{D_2}}\right)\right].
\]
For a fixed $\epsilon>0$ and each $x\in D_2$ we put
\[
r(x)=\frac12\rho_{D_2}(x)\wedge\epsilon \quad \mbox{ and }\quad
 B(x)=B(x, r(x)).
\]
Define a sequence of stopping times $\{T_m,
m\ge1\}$ as follows:
\[
T_1=\inf\{t>0: Y_t\notin B(Y_0)\},
\]
and for $m\ge 2$,
$$
T_m= T_{m-1}+
\tau_{B(Y_{T_{m-1}})}\circ\theta_{T_{m-1}}
$$
if $Y_{T_{m-1}}\in D_2$,  and $T_m= \tau^Y_{D_2}$ otherwise. Then
$\{h(Y_{T_m}), m\ge1\}$ is a martingale under $P_x$ for any $x\in
D_2$. Since $ D_2$ is smooth,  we know from (4.3) in \cite{KS} that
$\P_x(Y_{\tau^Y_{D_2}} \in \partial D_2)=0.$ Thus we have
$P_x(\tau^Y_{D_2}=T_m\, \mbox{ for some $m\ge1$})=1$. Since $h$ is
bounded on $D_2$, we have
\[
\left|E_x\left[h(Y_{T_m});\, T_m<\tau^Y_{D_2}\right]\right|\,\le\,
c\, P_x\left(T_m<\tau^Y_{D_2}\right)\,\rightarrow\, 0.
\]
Take a domain $D_3$ such that $\overline{D_2}\subset D_3
\subset\overline{D_3}\subset D$, then $h$ is continuous and
therefore bounded on $\overline{D_3}$. By Lemma \ref{l:MH1}, we have
$E_x[h(Y_{\tau^Y_{D_2}})] <\infty$. Thus by the dominated
convergence theorem
\[
\lim_{m\to\infty}E_x\left[h\left(Y_{\tau^Y_{D_2}}\right);\,
T_m=\tau^Y_{D_2}\right]=E_x\left[h\left(Y_{\tau^Y_{D_2}}\right)\right].
\]
Therefore
\begin{eqnarray*}
h(x)& =& \lim_{m\to\infty}E_x\left[h(Y_{T_m})\right] \\
&=& \lim_{m\to\infty}E_x\left[h(Y_{\tau^Y_{D_2}});\,
T_m=\tau^Y_{D_2}\right] + \lim_{m\to\infty} E_x\left[h(Y_{T_m});\,
T_m<\tau^Y_{D_2}\right]\, = \, E_x\left[h(Y_{\tau^Y_{D_2}})\right].
\end{eqnarray*}
\qed

Recall that a point $z\in \partial D$ is said to be a regular
boundary point for $Y$ if $P_z(\tau^Y_D=0)=1$ and an irregular
boundary point if $P_z(\tau^Y_D=0)=0$. It is well known that  if
$z\in \partial D$ is regular for $Y$, then for any $x\in D$,
$G^Y_D(x, y)\rightarrow 0$ as $y\rightarrow z$.

\begin{lemma}\label{lemma:boy0}
\begin{description}
\item{(1)} If
$z, w\in\partial D$, $z\neq w$ and $w$ is a regular boundary point for $Y$,
then $M^Y_D(x, z)\to 0$ as
$x\to w$.\item{(2)}
The mapping
$(x, z)\mapsto M^Y_D(x, z)$ is continuous on
$D\times\partial D$.
\end{description}
\end{lemma}

\pf Both of the assertions can be proved easily using our Lemma
\ref{l:5B} and Theorem \ref{t2.2}. We skip the proof since the
argument is almost identical to the one on  page 235 of \cite{B2}.
\qed

So far we have shown that the Martin boundary of $D$ can be
identified with a subset of the Euclidean boundary $\partial D$.
The main result of this section is as follows:

\begin{thm}\label{t3.1}
The Martin boundary and the minimal Martin boundary of $D$ with
respect to $Y$ can be identified with the Euclidean boundary of $D$.
\end{thm}

\pf For each fixed $z\in \partial D$ and $x\in D$, recall the Martin
kernel for $X^D$ from \cite{SW}:
\begin{equation}\label{e:t3.0}
M^X_D(x, z)\,:=\,\lim_{D\ni y \to z}\frac{G^X_D(x, y)}{G^X_D(x_0,
y)}.
\end{equation}
By Theorem \ref{t:ge}, there exists a constant $c >0$ such that
\begin{equation}\label{e:t3.1}
c^{-1} M^X_D(x, z) \,\le \, M^Y_D(x, z) \,\le \, c M^X_D(x, z),
\quad  (x,z)\in D \times \partial D.
\end{equation}
It is shown in pages 471-472 of \cite{SW} that if $M^X_D(\cdot, z_1)
=M^X_D(\cdot, z_2)$
then $\delta_{z_1}=\delta_{z_2}$.
In fact, by following the same proof, we can get
that if $c^{-1} M^X_D(\cdot, z_2)
\le M^X_D(\cdot, z_1)
\le cM^X_D(\cdot, z_2)$
then $c^{-1} \delta_{z_2} \le \delta_{z_1} \le c \delta_{z_2}$.
Thus by (\ref{e:t3.1}), if  $M^Y_D(\cdot, z_1)
=M^Y_D(\cdot, z_2)$, then $c^{-1} \delta_{z_2} \le
\delta_{z_1} \le c \delta_{z_2}$,
which implies that $z_1=z_2$.

We have shown that each Euclidean boundary point corresponds to a
different nonnegative harmonic function for $Y^D$, hence the Martin
boundary of $D$ for $Y$ coincides with the Euclidean boundary
$\partial D$. Finally we will show that, for every $z\in \partial
D$, $M^Y_D(\cdot, z)$ is a minimal harmonic function for $Y^D$,
hence the minimal Martin boundary of $D$ can be identified with the
Euclidean boundary.

Suppose that  $x\mapsto M^Y_D(x, z_0) $ is not a minimal harmonic
function of $Y^D$ for some $z_0\in \partial D$. Then there is a
non-trivial harmonic function $h\geq 0$ of $Y^D$ such that $h(x)
\leq M^Y_D(x, z_0)$ but $h$ is not a constant multiple of $M^Y_D(x,
z_0)$.

We know from
the general theory in Kunita and Watanabe \cite{KW} that non-negative harmonic
functions for $Y^D$ admit a Martin representation.
Since the Martin boundary of $D$ for $Y$ coincides with the
Euclidean boundary $\partial D$, there is a finite measure $\nu$ on
$\partial D$ which is not concentrated at $z_0$ such that
$$
h(x) = \int_{\partial D} M^Y_D(x,w) \nu (dw), \quad \mbox{ for }
x\in D.
$$
Define $u(x):= \int_{\partial D} M^X_D(x,w) \nu (dw) $, which is a
non-trivial positive harmonic function for $X^D$. By (\ref{e:t3.1})
and our assumption on $h$, $u(x) \leq c^2 \, M^X_D(x, z_0)$ for all
$x\in D$. Since $x\mapsto M^X_D(x, z_0)$ is a minimal harmonic
function for $X^D$ (see page 418 of \cite{SW}), $u$ has to be a
constant multiple of $M^X_D(\cdot, z_0)$. By the uniqueness in the
Martin representation for $X$ (see (4.1) in \cite{SW}), $\nu$ has to
be concentrated at point $z_0$. This is a contradiction. This proves
that $x\mapsto M^Y_D(x, z)$ is  a minimal harmonic function of
$Y^D$ for every $z\in
\partial D$. \qed

As a consequence of Theorem \ref{t3.1}, we conclude that
for every nonnegative harmonic function $h$ for $Y^D$, there exists a unique
finite measure $\mu$ on $\partial D$ such that
\begin{equation}\label{e3.1}
h(x)=\int_{\partial D}M^Y_D(x, z)\mu(dz), \qquad x\in D.
\end{equation}
We call $\mu$ the Martin measure of $h$.

By Theorem \ref{t:geC}, we get the following sharp estimates on
Martin Kernel.

\begin{thm}\label{t:MKE}
If $D$ is a bounded roughly connected $C^{1,1}$ open set, there
exists $c:=c(x_0, D, \alpha)>0$ such that
$$ \frac1{c} \frac {\rho_D (x)^{\alpha/2}}
{|x - z |^{d}} \,\leq\, M^Y_D (x, z) \,\leq\, c \,\frac {\rho_D
(x)^{\alpha/2}}{| x - z |^{d}} .
$$
\end{thm}

\section{Relative Fatou Type Theorem}

In  this section, we establish a relative Fatou type theorem for
truncated stable processes in bounded roughly connected $\kappa$-fat
open set. Throughout this section, we assume that $D$ is a bounded
roughly connected $\kappa$-fat open set in $\R^d$. The arguments of
this section are similar to those in \cite{K2}. We spell out some of
the details for the readers' convenience.

In the previous section, we have shown that the Martin kernel
$M^Y_D(x , z)$ is harmonic for $Y^D$. In fact, a stronger result is
true.

\begin{lemma}\label{lemma:boy}
For each $z \in \partial D$,
 $M^Y_D(\, \cdot \, , z)$ is a bounded regular harmonic function of $Y$
in $D \setminus B(z, \eps)$ for every $\eps > 0$.
\end{lemma}

\pf Fix $z \in \partial D$ and $\eps > 0$, and let $h(x):=M^Y_D(x,
z)$ for $x \in \RR^d$. By (\ref{e:t3.1}),  $h$ is bounded on  $\RR^d
\setminus B(z, \eps/2)$. Take an increasing sequence of smooth open
sets $\{ D_m \}_{m \ge 1}$ such that $\overline{D_m} \subset
D_{m+1}$ and $\cup^{\infty}_{m=1} D_m = D \setminus B(z, \eps)$. Set
$\tau_m := \tau^Y_{D_m}$ and $\tau_{\infty} := \tau^Y_{D \setminus
B(z, \eps)}$ . Then $\tau_m \uparrow \tau_{\infty}$. Since $Y^D$ is
a Hunt process, $\lim_{m \to \infty} Y^D_{\tau_m} =
Y^D_{\tau_\infty} $ by the quasi-left continuity of $Y^D$. Set $A =
\{ \,\tau_m = \tau_\infty  ~\mbox{ for some } m \ge 1 \}$. Let $N$
be the set of irregular boundary points of $D$ for $Y$. $N$ is
semi-polar by Proposition II.3.3 in \cite{BG}, which is polar in our
case (Theorem 4.1.2 in \cite{FOT}).
 Thus
\begin{equation}\label{zero1}
\P_x(Y_{\tau_\infty} \in N)=0,  ~~~  x \in D.
\end{equation}
By Lemma \ref{lemma:boy0},
 if $w \in \partial D, w \not= z $
and $w$ is a regular boundary point, then
$
h(x) \to 0 $ as  $ x \to w$
so that $h$ is continuous on $\overline{D \setminus B(z, \eps)} \setminus N$.
Therefore, since $h$ is bounded on $\RR^d \setminus B(z, \eps/2)$,
by the bounded convergence theorem and (\ref{zero1}), we have
\begin{eqnarray*}
&&\lim_{m \to \infty} \E_x \left[ \,h (Y^D_{\tau_m}) \,;\, \tau_m <
\tau_\infty  \right] =\lim_{m \to \infty} \E_x \left[ \,h
(Y_{\tau_m})1_{ \overline{D \setminus B(z, \eps)} \setminus
N}(Y_{\tau_\infty})
\,;\, \tau_m < \tau_\infty \right]\\
&&= \E_x \left[ \,h (Y_{\tau_\infty})1_{\overline{D \setminus B(z,
\eps)} \setminus N}(Y_{\tau_\infty}) \,;\, A^c\,\right] = \E_x
\left[ \,h (Y_{\tau_\infty}) \,;\, A^c\,\right].
\end{eqnarray*}
By the boundedness of $h$ on  $\RR^d \setminus B(z, \eps/2)$, we can
find two smooth open sets $U_1$ and $U_2$ such that $\overline{D
\setminus B(z, \eps)} \subset U_1 \subset \overline{U}_1 \subset
U_2$ and $h$ is the bounded on $\overline{U}_2$. Pick a point $z_0
\in D \setminus U_2$. Since
$$
|v - w| \,\ge\, \frac{\mbox{dist}(U_1, D\setminus
U_2)}{\mbox{diam}(D)} |v - z_0| \,\ge\, \left(\frac{\mbox{dist}(U_1,
D\setminus U_2)}{\mbox{diam}(D)}\right)^2 |v - w|,
 \qquad \forall (v,w) \in (U_1 \cap D)\times (D \setminus U_2),
$$
we have, by (\ref{t_Levy})-(\ref{e:pkY}),
\begin{eqnarray*}
\E_x \left[\,h \left(Y_{\tau^Y_{U_1 \cap D}}\right) 1_{ \{
Y_{\tau^Y_{U_1\cap D}} \in D \setminus U_2 \}} \right] &=& {\cal
A}(d, -\alpha)\int_{D \setminus U_2} h(z) \int_{U_1 \cap D \cap \{
|v-z| <1 \}}
\frac{G^Y_{U_1 \cap D} (x, v)}{|v-z|^{d + \alpha}} dvdz \\
&\le& c_1\int_{D \setminus U_2} h(z)dz \int_{U_1 \cap D \cap \{
|v-z| <1 \}} \frac{G^Y_{U_1 \cap D} (x, v)}
{|v-z_0|^{d + \alpha}} dv\\
&\le& c_2  \,||h||_{L^1 (D)} \int_{D \setminus U_1} \int_{U_1 \cap
D\cap \{ |v-z| <1 \}}   \frac{G^Y_{U_1
\cap D} (x, v)}{|v-z|^{d + \alpha}} dvdz\\
&=& c_2 \,||h||_{L^1 (D)} \P_x \left( Y_{\tau^Y_{U_1\cap D}} \in D
\setminus U_1 \right) \,<\, \infty .
\end{eqnarray*}
Therefore
$$
\E_x \left[ \,h ( Y_{\tau_\infty} );\, \tau_\infty < \tau^Y_D
\,\right] \, \le\, \E_x \left[ \,h (Y_{\tau_\infty}) 1_{ \{
Y_{\tau_\infty} \in U_2 \cap D\} } \,\right] \, +\, \E_x \left[\,
h\left(Y_{\tau^Y_{U_1\cap D}}\right) 1_{ \{ Y_{\tau^Y_{U_1\cap D}}
\in D \setminus U_2 \} } \,\right] < \infty .
$$
Thus by the dominated convergence theorem,
\begin{eqnarray*}
h(x) &=& \lim_{m \to \infty} \E_x \left[ \,h (Y_{\tau_m});
\,\tau_m< \tau^Y_D\right]\\
&=&  \lim_{m \to \infty} \E_x \left[\,h (Y_{\tau_\infty}) ; \,\tau_m
= \tau_\infty < \tau^Y_D \right]\, + \,\lim_{m \to \infty} \E_x
\left[ \,h (Y_{\tau_m}) ;
\,\tau_m < \tau^Y_\infty< \tau^Y_D\right]\\
&=& \E_x \left[\, h (Y_{\tau_\infty});\,\tau^Y_\infty< \tau^Y_D \right] .
\end{eqnarray*}
\qed

\begin{prop}\label{prop:3.2}
For every $  \lambda \in (0,1)$, there exists $c = c ( D, \alpha,
\lambda)>0$ such that if $y \in D $ and $|y - x_0 | > 2 \rho_D (y) $
then
\begin{equation}\label{eqn:3.5}
 \P_{x_0} \left(T_{B^{\lambda}_y}  < \tau^Y_D \right)
\,\ge\, c \,G^Y_D(x_0, y) \rho_D (y)^{d- \alpha},
\end{equation}
where $B^{\lambda}_y := B(y, \lambda \rho_D (y))$
and $T_{B^{\lambda}_y} = \inf\{t>0: \, Y_t \in B^{\lambda}_y \}$.
\end{prop}

\pf First note that $x_0 \not\in B(y, \rho_D(y))$. Since $G^X_D(x_0
, \,\cdot\,)$ is harmonic for $X$ in $D \setminus \{x_0\}$, by the
Harnack inequality for $X$ (Lemma 2 in \cite{B}) and Theorem
\ref{t:ge}, there exist $c_i = c_i ( D, \alpha, \lambda)>0$, $i=1,
2, 3$, such that
\begin{equation}\label{eqn:3.6}
 \int_{B^{\lambda}_y} G^Y_D(x_0 , z) dz\,\ge\,
c_1\int_{B^{\lambda}_y} G^X_D(x_0 , z)  dz\,\ge\, c_2 \, G^X_D(x_0 ,
y) \rho_D (y)^d\,\ge\, c_3 \, G^Y_D(x_0 , y) \rho_D (y)^d .
\end{equation}
Using the strong Markov property, one can easily see that
\begin{equation}\label{eqn:3.7}
 \int_{B^{\lambda}_y} G^Y_D(x_0 , z) dz \,\le\,
  \left(\sup_{w \in
\overline{B^{\lambda}_y}} \E_w \int^{\tau^Y_D}_0 1_{B^{\lambda}_y}
(Y_s) ds\right)\P_{x_0}\left(T_{B^{\lambda}_y} <\tau^Y_D \right) .
\end{equation}
On the other hand, by (\ref{e:G_ub}),
\begin{equation}\label{eqn:3.8}
\E_w \int^{\tau^Y_D}_0 1_{B^{\lambda}_y} (Y_s) ds
\,=\,\int_{B^{\lambda}_y} G_D^Y(w,v)dv
\,\le\,
c_4\int_{B^{\lambda}_y} \frac{dv}{|w-v|^{d- \alpha}}\,\le\, c_5
\,\rho_D (y)^{\alpha}
\end{equation}
for every $w \in \overline{{B^{\lambda}_y}}$ for
some $c_i=c_i(D, \alpha, \lambda)>0$, $i=4, 5$.
Combining (\ref{eqn:3.6})-(\ref{eqn:3.8}),
we have (\ref{eqn:3.5}).
\qed

\begin{defn}\label{D:5.1}
A nonnegative Borel measurable function
$f$ defined on $D$
is said to be

\begin{description}
\item{(1)} excessive with respect to $Y^D$ if for every $x\in D$ and $t>0$,
$$
\E_x\left[  f(Y^D_t)\right] \,\le\, f(x)\quad \mbox{and}\quad
\lim_{t \downarrow 0} \E_x\left[  f(Y^D_t)\right] \,=\, f(x);
$$
\item{(2)} superharmonic with respect to $Y^D$ if
$f$ is lower semi-continuous in $D$, and $$ f(x)\,\ge \,\E_x\left[
f\left(Y^D_{\tau^Y_{B}}\right)\right], \qquad x\in B,
$$
for every open set $B$ whose closure is a compact subset
of $D$.
\end{description}
\end{defn}

Suppose $h>0$ is a positive superharmonic function for $Y^D$. Since
$Y^D$ is a Hunt process satisfying the strong Feller property (i.e.,
for every $f \in L^{\infty}(D)$, $\E_x[f(Y^D_t)]$ is bounded and
continuous in $D$), $h$ is excessive for  $Y^D$ (for example,  see
\cite{CK2}). For any positive superharmonic function $h$ for $Y^D$,
let $D_h:=\{x \in D:\,h(x) < \infty \}$ and define
$$
p^h_D(t,x,y)\,:=\, h(x)^{-1} p_D^Y(t,x,y) h(y), \quad t>0, ~x,y \in D_h.
$$
Then $p^h_D(t,x,y)$ is a transition density and it determines a nice
Markov process on $D_h \cup \{ \partial \}$ (for example, see
\cite{KW}). This process is called an $h$-conditioned truncated
stable process and we will use $\E^h_x$ to denote the expectation
with respect to this process.

Let $\{{\cal F}_t, t \ge 0 \}$  be the minimal admissible
$\sigma$-fields generated by $Y$.
For any stopping time
$T$ of $\{{\cal F}_t , t \ge0\}$, ${\cal F}_{T+}$ is the class of subsets $A$ of ${\cal F}$
such that $A \cap \{T \le t\}  \in {\cal F}_{t+}$ for every $t>0.$
Similar to Propositions 5.2-5.4 of
\cite{CZ} (also see Lemmas 3.11-12 in \cite{CS3}), we have the
following.
\begin{lemma}\label{t:ht}
For any stopping time $T$ and any ${\cal F}_{T+}$-measurable function
$\Psi \ge 0$,
$$
\E_x^h\left[\Psi; T< \tau^Y_D \right]\,=\, h(x)^{-1} \E_{x}
\left[\Psi \cdot h(Y_T); T< \tau^Y_D\right], \quad x \in D.
$$
\end{lemma}

Let $(\P^z_x , Y^D_t)$ be  the $h$-conditioned truncated stable
process with $h(\cdot)=M^Y_D(\cdot , z)$. The following theorem is
known for symmetric stable processes (see \cite{K2}).

\begin{thm}\label{thm:son1}
$$
\P^z_x \left( \lim_{t \uparrow \tau^Y_D} Y^D_t = z , \, \tau^Y_D<
\infty \right) = 1, ~~~\mbox{ for every } x \in D,~ z \in \partial
D.
$$
\end{thm}

\pf
By Theorem 3.15 in \cite{KS6}, we know that
$$
\sup_{(x,z) \in D \times \partial D}\E^z_x[\tau^Y_D] \,< \,\infty.
$$
Therefore
$
\P^z_x ( \tau^Y_D< \infty ) = 1$ for every  $x \in D$ and
$z \in \partial D.$

Now we fix $x \in D$ and $z \in \partial D$. We claim that $\P^z_x (
\lim_{t \uparrow \tau^Y_D } Y^D_t = z ) = 1$. Let $r_m= 1/2^m$,
$B_m:=B(z,r_m)$, $D_m:= D \setminus \overline{B_m}$ and set $T_m :=
\inf\{t>0; Y^D_t \in B_m\}$ and $R_m = \tau^Y_{B_m \cap D}$. We may
suppose that $x \in D_m$. Since $\P_x(Y_{\tau^Y_{D_m}} \in
\partial B_m)=0$ by (4.4) in \cite{KS}, by Lemmas \ref{lemma:boy} and \ref{t:ht}, we have
$$
M^Y_D(x,z) \,=\,\E_x \left[ M^Y_D ( Y_{\tau^Y_{D_m}}, z ) \right]
\,=\,\E_x \left[ M^Y_D ( Y_{T_m}, z ) ;\, T_m < \tau^Y_D \right]
\,=\, M^Y_D(x,z)\, \P^z_x (T_m < \tau^Y_D) .
$$
It follows that for all $m \ge 1$ we have $\P^z_x (T_m < \tau^Y_D) = 1$.
Let $L_k := \sup_{y \in B^c_k \cap D} M^Y_D(y,z)$, which
is finite by  Lemma \ref{lemma:boy}.
For $k < m$ (see \cite{CS3}, Theorem 3.17),
$$
\P^z_x \left[ T_m < \tau^Y_D ,~ R_k \circ \theta_{T_m} < \tau^Y_D
\right] \,\le\, \frac{L_k}{M^Y_D(x,z)}\P_x\left(T_m <
\tau^Y_D\right)
$$
Using Theorem \ref{t:ge} one can easily show that
any singleton $\{z\}$ in $\R^d$ has zero capacity with respect to $Y$.
Thus we have
$$
\limsup_{m \to \infty} \P_x \left(T_m < \tau^Y_D\right) \, \le\,
\P_x \left(T^Y_{\{z\}} \le \tau^Y_D \right) \,\le\, \P_x
\left(T^Y_{\{z\}} < \infty \right) \,=\, 0,
$$
where $T^Y_{\{z\}} := \inf\{t>0: Y_t =z\}$. The rest of the proof is
similar to the corresponding part of the proof in Theorem 3.17 in
\cite{CS3}. We skip the details. \qed

The theorem above implies that for every Borel subset $K \subset
\partial D$,
$$
\P^{z}_x \left( \lim_{t \uparrow \tau^Y_D} Y^D_t \in K \right) =
1_K(z),\qquad \forall (x,z) \in D \times \partial D.
$$
So the next theorem follows easily from the Martin representation
for $Y^D$ in (\ref{e3.1}).

\begin{thm}\label{thm:son2}
Let $h$ be a positive harmonic function for $Y^D$ with the Martin
measure $\nu$. Then for any $ x \in D$ and any Borel subset $K$ of
$\partial D$, we have
$$
\P^h_x \left(\lim_{t \uparrow \tau^Y_D} Y^D_t \in K\right) =
\frac1{h(x)} \int_K M^Y_D(x,w) \nu(dw).
$$
\end{thm}

The proof of the next proposition is identical to the proof of Proposition 3.5 in \cite{K2}. So we skip
the proof.

\begin{prop}\label{prop:3.3}
Let $h$ be a positive harmonic function for $Y^D$ with the Martin
measure $\nu$ satisfying $\nu(\partial D)=1$. If $A \in {\cal
F}_{\tau^Y_D}$, then for every Borel subset $K$ of $\partial D$,
$$
\int_K \P^z_{x_0} (A)\nu(dz) \,=\, \P^h_{x_0}\left(A \cap
\left\{\lim_{t \uparrow \tau^Y_D} Y^D_t \in K\right\} \right).
$$
\end{prop}

\begin{defn}\label{def:3.4}
$A \in {\cal F}_{\tau^Y_D}$ is said to be shift-invariant if whenever
$T < \tau^Y_D $ is a stopping time, $ 1_A \circ \theta_T = 1_A$
$\P_x$-a.s. for every $x \in D$.
\end{defn}

\begin{prop}[0-1 law]\label{prop:3.5}
 If $A$ is shift-invariant, then $x \to \P^z_x (A)$ is a
constant function which is either $0$ or $1$.
\end{prop}
\pf
See the proof of Proposition 3.7 in \cite{K2}. \qed

We recall the definition of Stolz open sets for a $\kappa$-fat open
set from \cite{K2}. Recall that $(R,\kappa)$ is the characteristic
of the $\kappa$-fat open set $D$.

\begin{defn}\label{Stolz}
For $z \in \partial D$ and $\beta > (1-\kappa)/\kappa$, let
$$
A^{\beta}_z := \left\{ y \in D ;~ \rho_D(y) \,<\,
\frac{\rho_D (x_0)}{3}\wedge R ~\mbox{ and } ~|y-z|
<  \beta \rho_D(y) \right\}.
$$
We call $A^{\beta}_z$ the Stolz open set for $D$ at $z$
with the angle $\beta > (1-\kappa)/\kappa$.
\end{defn}

We know from Lemma 3.9 of \cite{K2} that for every $z
\in \partial D$ and $\beta > (1-\kappa)/\kappa$, there
exists a sequence $\{y_k\}_{k \ge 1} \subset A^{\beta}_z$
such that $\lim_{k \to \infty} y_k =z$.

\begin{prop}\label{prop:3.6}
Given $z \in \partial D$, $\lambda \in (0,1)$ and $\beta > (1-\kappa)/\kappa$,
there exists $c = c(D, \alpha , \lambda, x_0, \beta) > 0$
such that if $y \in A^{\beta}_z$
then
$$
\P^z_{x_0} \left(T_{B^{\lambda}_y} < \tau^Y_D \right) ~>~ c,
$$
where $B^{\lambda}_y := B(y, \lambda \rho_D (y))$
and $T_{B^{\lambda}_y} := \inf\{t>0: \, Y_t \in B^{\lambda}_y \}$.
\end{prop}

\pf Fix $z \in \partial D$ and  $\beta > (1-\kappa)/\kappa$. Recall
 from (\ref{e:t3.0}) that $M^X_D (x,z)$ is the Martin kernel for
$X^D$. Since $M^X_D(\cdot,z)$ is harmonic for $X^D$ in $D$, by the Harnack
inequality for $X$ (Lemma 2 in \cite{B}) and Proposition
\ref{prop:3.2} we have
$$
\E_{x_0} \left[ M^X_D\left(Y_{T_{B^{\lambda}_y}} , z\right) ;
 T_{B^{\lambda}_y} < \tau^Y_D \right]
\geq c_1\P_{x_0} \left(  T_{B^{\lambda}_y} < \tau^Y_D \right)\,
M^X_D(y,z)  \geq c_2  G^Y_D(x_0, y) \rho_D (y)^{d-\alpha}
M^X_D(y,z).
$$
Thus by Lemma \ref{t:ht}
 and (\ref{e:t3.1}),
\begin{eqnarray*}
&&\P^z_{x_0} \left( T_{B^{\lambda}_y} < \tau^Y_D \right) = \E_{x_0}
\left[ M^Y_D\left(Y_{T_{B^{\lambda}_y}} , z\right) ;\,
 T_{B^{\lambda}_y} < \tau^Y_D \right]\\
&&\geq c_3\,\E_{x_0} \left[ M^X_D\left(Y_{T_{B^{\lambda}_y}} ,
z\right) ;\,
 T_{B^{\lambda}_y} < \tau^Y_D \right] \,\ge\,
 c_4
\, G^Y_D(x_0, y) \,\rho_D (y)^{d-\alpha}\, \lim_{w \in D \rightarrow
z} \frac{G^X_D(y,w)}{G^X_D(x_0,w)}.
\end{eqnarray*}
By Theorem \ref{t:ge}, the last quantity  above is
greater than or equal to
$$
c_5 \,G^X_D(x_0, y) \, \rho_D (y)^{d-\alpha}\, \lim_{D \ni w \to z}
\frac{G^X_D(y,w)}{G^X_D(x_0,w)},
$$
which is bounded below by a positive constant by the last
display in the proof of
Proposition 3.10 in \cite{K2}.
\qed

To establish our relative Fatou type theorem, we need to control
the oscillation of nonnegative harmonic functions of $Y$. The
following three lemmas are preparations for our
oscillation result, Proposition \ref{prop:3.7} below.

We first recall from \cite{KS} that if $U$ is an open set with $\mbox{diam}(U) \le \frac12$,
the Feynman-Kac semigroup $(Q^U_t)$ defined by
$$
Q^{U}_tf(x):= \E_x\left[\exp(\int^t_0q_U(X^{U}_s)ds)f(X^{U}_t)\right]
$$
is the semigroup of $Y^U$, where
$$
q_U(x):= \int_{U^c\cap\{|x-y|\ge 1\}}\nu^X(x-y)dy={\cal A} (d, - \alpha)
   \int_{U^c\cap\{|x-y|\ge 1\}} |x-y|^{-(d+\alpha)} dy.
$$
Note that
\begin{equation}\label{cal}
0\le q_U(x) \le {\cal A} (d, - \alpha)
\int_{\{|x-y|\ge 1\}} |x-y|^{-(d+\alpha)}dy=:{\cal B}(d, \alpha),
\quad \forall x\in D.
\end{equation}

The next lemma is a modification of Proposition 3.2 in \cite{KS}.
Recall that $r_0 \in (0, \frac14)$ is the constant from Theorem \ref{T:Har0}.

\begin{lemma}\label{G_1}
For every $\eps >0$, there exists a positive constant
$r_1=r_1(\alpha, d, \eps) \le r_0$ such that for all
$r\in (0, r_1]$ and $a \in \R^d$, we have
$$
G^Y_{B(a, r)}(x, y)\,\le\, \sqrt{1+\frac{\eps}{4}}\, G^X_{B(a,
r)}(x, y), \quad x, y\in B(a, r).
$$
\end{lemma}

\pf Let $B_r:=B(0,r)$ with $r \le \frac14$. For any $z \in B_r$, let
$(\P^z_x , X^{B_r}_t)$ be the $G_{B_r}(\cdot , z)$-transform of
$(\P_x , X^{B_r}_t)$, that is, for any nonnegative Borel functions
$f$ in $B_r$,
$$\E^z_x \left[f(X^{B_r}_t)\right] = \E_x \left[\,
\frac{G^X_{B_r}(X^{B_r}_t , z)}
{G^X_{B_r}(x,z)}f(X^{B_r}_t)\right].
$$
It is well known that there exists a positive constant $C$
independent of $r$ such that
\begin{equation}\label{3g4balls}
\frac{G^X_{B_r}(x, y)G^X_{B_r}(y, z)}
{G^X_{B_r}(x, z)}\,\le\, C\,(|x-y|^{\alpha-d}+|y-z|^{\alpha-d}),
\quad \forall\, x, y, z\in B_r.
\end{equation}
So there exists a positive constant $r_1$ such that for any $r\in (0, r_1]$
and all $x, z\in B_r$,
\begin{equation}\label{e:3G_ball}
{\cal B}(d, \alpha)\,\E^z_x\tau^X_{B_r}\, = \,{\cal B}(d,
\alpha)\int_{B_r} \frac{G^X_{B_r}(x, y)G^X_{B_r}(y, z)} {G^X_{B_r}(x,
z)}dy \,<1- \frac1{\sqrt{1+\eps/4}},
\end{equation}
 where ${\cal B}(d, \alpha)$ is the constant in (\ref{cal}).
Hence by (\ref{cal}) and Khasminskii's lemma (see, for instance,
Lemma 3.7 in \cite{CZ})
we get that for
$r \in (0, r_1]$
$$
\E^z_x\left[\exp(\int^{\tau^X_{B_r}}_0q(X^{B_r}_s)ds)\right]\,\le\,
\E^z_x  \left[ \exp\left({\cal B}(d, \alpha)
\tau^X_{B_r}\right)\right]\,\le\, \sqrt{1+\frac{\eps}{4}}.
$$
Since
$$
G^Y_{B_r}(x, z)\,=\,G^X_{B_r}(x, z)\,\E^z_x\left[
\exp(\int^{\tau^X_{B_r}}_0q(X^{B_r}_s)ds)\right], \quad x, z\in B(0,
r),
$$
using the translation invariance property of our Green functions, we
arrive at our desired result. \qed

The above lemma implies that
$$
K_{B(a,r)}^Y(x,z)\,\le\,
{\cal A}(d, -\alpha)\int_{B(a,r)}
 \frac{G_{B(a,r)}^Y(x,y)}{|y-z|^{d+\alpha}}
 dy\,\le\, {\cal A}(d, -\alpha)\sqrt{1+\frac{\eps}4}\int_{B(a,r)}\,
 \frac{G_{B(a,r)}^X(x,y)}{|y-z|^{d+\alpha}}
 dy
$$
for every $z \in A(a, r, \infty)$ and $x \in B(a,r)$.
Thus we have proved following lemma.

\begin{lemma}\label{l:P_0}
For every $\eps >0$, there exists a positive constant
$r_1=r_1(\alpha, d, \eps) \le r_0$ such that
for $r\in (0, r_1]$  and $z \in A(a, r, \infty)$,
\begin{equation}\label{CC3}
K_{B(a,r)}^Y(x,z) \,\le\, \sqrt{1+\frac{\eps}4}\,
K^X_{B(a,r)}(x,z),\qquad x \in B(a,r).
\end{equation}
\end{lemma}

\begin{lemma}\label{l:P_1}
For every $\eps >0$, there exist  positive constants
$r_1=r_1(\alpha, d, \eps) \le r_0$  and $\lambda_1 =\lambda_1 (\alpha, d, \eps)<\frac12$ such that
for every  $r < r_1$,  $z \in A(a, r, 1-r)$
and $x_1, x_2 \in B(a,\lambda_1 r )$,
$$
 \frac1{1+{\eps}/2} \,K^Y_{B(a,r)}(x_2,z) \,\le\, K_{B(a,r)}^Y(x_1,z)
\,\le\,  (1+\frac{\eps}2)\, K^Y_{B(a,r)}(x_2,z).
$$
\end{lemma}

\pf
By (4.7) in \cite{KS} and Lemma \ref{l:P_0} above, there exist
$r_1=r_1(\alpha, d, \eps) \le r_0 <1/4$ such
that for every $r\le r_1$ and $z \in A(a, r, 1-r)$,
\begin{equation}\label{aaa1}
K^X_{B(a,r)}(x,z) \,\le\, K_{B(a,r)}^Y(x,z) \,\le
\,\sqrt{1+\frac{\eps}4}\, K^X_{B(a,r)}(x,z), \quad x \in B(a,r).
\end{equation}
By the explicit formula for  $K^X_{B(a,r)}(x,z)$ in (\ref{P_f}), we
see that for $x_1, x_2
\in B(a, \lambda r)$ and $z \in A(a, r, 1-r)$
\begin{equation}\label{aaa2}
K^X_{B(a,r)}(x_1,z) \le  \frac{(1+\lambda)^d}{(1-\lambda)^d
(1-\lambda^2)^{\alpha/2}}
K^X_{B(a,r)}(x_2,z).
\end{equation}
The lemma follows easily from the inequalities (\ref{aaa1})-(\ref{aaa2}).
\qed

Now we can prove our oscillation result for nonnegative harmonic functions
of $Y$.

\begin{prop}\label{prop:3.7}
For any given $\eps >0$, there exists $ \lambda_0 =
\lambda_0 (\eps, \alpha, d) \in (0,1/4)$
such that whenever $u$ is a positive harmonic
function for $Y$ in $D$,
$$
\frac{1}{1+ \eps }u(z)\, \le\, u(a)\, \le\, (1+ \eps ) u(z)
$$
for every $a \in D$ with $\rho_D (a)< r_1$
and $ z \in B^{\lambda_0}_a:= B(a, \lambda_0 \rho_D (a))$, where
$r_1$ is the constant from Lemma \ref{l:P_1}.
\end{prop}
\pf
Recall that $\lambda_1 <\frac{1}{2}$ is the
constant from Lemma \ref{l:P_1}.
Fix a positive harmonic function  $u$ for $Y$ in $D$ and a point
$a \in D$ with $\rho_D (a)< r_1$.
For any $\eta\in (0, 1)$, we put $B_a^{\eta}=
B(a, \eta \rho_D (a))$. From now on we assume that
$\lambda \in (0,\frac14)$ and
put $\tau(\lambda) := \tau_{B^{\lambda}_a}$
and $\rho:=\rho_D (a)$.
By Theorem 4.1 in \cite{KS}, we have for any $y\in B^{\lambda \lambda_1}_a$,
\begin{eqnarray*}
u(y) \,=\, \E_{y} \left[u(Y_{\tau(\lambda)})\right]&=& \left(\int_{A(a, \lambda \rho, 1-\lambda \rho)} + \int_{A(a, 1-\lambda \rho, 1+ \lambda \rho/2)}\right)
K_{B^{\lambda}_a}^Y(y,z)u(z) dz \\
&&+
\E_{y} \left[u(Y_{\tau(\lambda)});\, Y_{\tau(\lambda)}
\in   A(a, 1+ \lambda \rho/2, 1+\lambda \rho)    \right].
\end{eqnarray*}
By Lemma \ref{l:P_1}, for $y, w \in B^{\lambda \lambda_1}_a$ we have
\begin{eqnarray}
\int_{A(a, \lambda \rho, 1-\lambda \rho)}
K_{B^\lambda_a}^Y(y,z) u(z)dz &\le& (1+\frac{\eps}2) \int_{A(a, \lambda \rho,
1-{\lambda \rho})} K_{B^\lambda_a}^Y(w,z) u(z)dz \nonumber\\
&=&(1+\frac{\eps}2)\, \E_{w} \left[u(Y_{\tau_{B(a,
\lambda \rho)}});\, Y_{\tau_{B(a, \lambda \rho)}}
\in   A(a, \lambda \rho, 1-\lambda \rho)    \right].\label{r_1}
\end{eqnarray}
Note that by Theorem 4.1 in \cite{KS},
$$
\P_{y} \left( Y_{\tau(\lambda)}
\in   A(a, 1+ \lambda \rho/2, 1+\lambda \rho)
, \, \tau(\lambda) =
\tau(\lambda/2)\right)  = \P_{y} \left( Y_{\tau(\lambda/2)}
\in   A(a, 1+\lambda \rho/2, 1+\lambda \rho)\right)
\, =\,0.$$
Thus by the strong Markov property, we have
\begin{eqnarray*}
&&\E_{y} \left[\,u(Y_{\tau(\lambda)});\, Y_{\tau(\lambda)}
\in   A(a, 1+\lambda \rho/2, 1+\lambda \rho)     \right]\\
&&=\E_{y} \left[\,u(Y_{\tau(\lambda)});\, Y_{\tau(\lambda)}
\in   A(a, 1+ \lambda \rho/2, 1+\lambda \rho) ,\,
\tau(\lambda) >
\tau(\lambda/2)   \right]\\
&&=
\E_{y} \left[\,    \E_{Y_{\tau(\lambda/2)}} \left[
u(Y_{\tau(\lambda)});\, Y_{\tau(\lambda)}
\in   A(a, 1+ \lambda \rho/2, 1+\lambda \rho)   \right]
 1_{ A(a,  \lambda \rho/2, \lambda \rho)}
\left(Y_{\tau(\lambda/2)}\right)
 \right].
\end{eqnarray*}
Let
$$
g(z):=\E_{z} \left[    u(Y_{\tau(\lambda)});\, Y_{\tau(\lambda)}
\in   A(a, 1+ \lambda \rho/2, 1+\lambda \rho)    \right]
$$
for $z \in  A(a, \lambda \rho/2, \lambda \rho)$,
and zero otherwise.
Then we have from the above argument that
$$
\E_{y} \left[u(Y_{\tau(\lambda)});\, Y_{\tau(\lambda)}
\in   A(a, 1+{\lambda \rho}/2, 1+\lambda \rho)     \right]=
\E_{y} \left[g(Y_{\tau(\lambda/2)})
 \right].
$$
Since the function
$y\mapsto \E_{y} [g(Y_{\tau(\lambda/2)})]$
is regular harmonic on $B(a, \lambda \rho/2)$ with respect to $Y$,
and is zero on $\overline{B(a, \lambda \rho)}^c$, we get
by Lemma \ref{l:P_1} and the argument in (\ref{r_1})  that for
$y, w \in B^{\lambda \lambda_1}_a$,
\begin{eqnarray}
&&\E_{y} \left[u(Y_{\tau(\lambda)});\, Y_{\tau(\lambda)}
\in   A(a, 1+{\lambda \rho}/2, 1+\lambda \rho)     \right]
\nonumber \\
&&\le \,(1+\frac{\eps}2) \,\E_{w} \left[
\E_{Y_{\tau(\lambda/2)}} \left[
u(Y_{\tau(\lambda)});\, Y_{\tau(\lambda)}
\in   A(a, 1+ \lambda \rho/2, 1+\lambda \rho)   \right]
 1_{ A(a,  \lambda \rho/2, \lambda \rho)}
\left(Y_{\tau(\lambda/2)}\right)
 \right] \nonumber \\
&&=\,(1+\frac{\eps}2)\,\E_{w} \left[u(Y_{\tau(\lambda)});
\, Y_{\tau(\lambda)}
\in   A(a, 1+{\lambda \rho}/2, 1+\lambda \rho)
 \right].\label{r_2}
\end{eqnarray}

On the other hand, by Lemma 4.5 in \cite{KS}, there exists $M=M(d, \alpha)>1$
such that
\begin{eqnarray*}
\int_{A(a, 1-\lambda \rho, 1+ \lambda \rho/2)}
K_{B^{\lambda}_a}^Y(y,z)u(z) dz
&\le& M \,\lambda^{\alpha}\,\rho^{\alpha}\int_{A(a,
1-\lambda \rho, 1+ \lambda \rho/2)}u(z) dz\\
&\le & M \, \lambda^{\alpha}\rho^{\alpha}\int_{A(a,
1- \rho/2, 1+ \rho/4)}u(z) dz.
\end{eqnarray*}
By applying Lemma 4.5 in \cite{KS} to
$K_{B^{1/2}_a}^Y(w,z)$ with
$w \in B^{\lambda \lambda_1}_a$ and $z \in A(a, 1- \rho/2,
1+\rho/4)$, we get
$$
\int_{A(a, 1-\rho/2, 1+\rho/4)}
\rho^{\alpha} u(z) dz \, \le\, M \,2^{\alpha}\int_{
A(a, 1-\rho/2, 1+ \rho/4)}
K_{B^{1/2}_a}^Y(w,z) u(z) dz,
$$
which is less than or equal to $M 2^{\alpha} u(w)$ by the harmonicity of $u$.
 Therefore we have
 \begin{equation}\label{r_3}
 \int_{A(a, 1-\lambda \rho, 1+ \lambda \rho/2)}
K_{B^{\lambda}_a}^Y(y,z)u(z) dz\, \le \, M^2 2^{\alpha} \lambda^{\alpha} u(w).
\end{equation}

Combining (\ref{r_1})-(\ref{r_3}) with
$\lambda:=(\eps^{\frac{1}{\alpha}} M^{-\frac{2}{\alpha}}
2^{-1-\frac{1}{\alpha}})\wedge\frac14$, we conclude that
for $y, w \in B^{\lambda \lambda_1}_a$
\begin{eqnarray*}
u(y)
&\le& (1+\frac{\eps}{2})
\E_{w} \left[u(Y_{\tau(\lambda)});\, Y_{\tau(\lambda)}
\in A(a, \lambda \rho, 1-\lambda \rho)
\cup  A(a, 1+\lambda \rho/2, 1+\lambda \rho)    \right]\,+ \,\frac{\eps}{2} u(w)\\
 & \le & (1+\eps) u(w).
\end{eqnarray*}
In particular,
$$
\frac{1}{1+ \eps } u(z)~\le~ u(a) ~\le ~(1+ \eps ) u(z)
$$
for every $ z \in B^{\lambda \lambda_1}_a$. The proposition is proved with
$ \lambda_0:=\lambda_1\left((\eps^{\frac{1}{\alpha}} M^{-\frac{2}{\alpha}}
2^{-1-\frac{1}{\alpha}})\wedge\frac14\right)$.
\qed

Now we are ready to establish a  relative Fatou type theorem of
harmonic function for $Y^D$. With Propositions \ref{prop:3.5} and
\ref{prop:3.7} in hand, the proof of the relative Fatou type theorem
is an easy modification of the proof of Theorem 3.13 in \cite{K2}.
We spell out detail for the readers' convenience.

\begin{thm}\label{T:Fatou}
Let $h$ be a positive harmonic function for $Y^D$ with the Martin
measure $\nu$.  If $u$ is a nonnegative harmonic function for $Y$ in
$D$, then for $\nu$-a.e. $z \in
\partial D$,
\begin{equation}\label{eqn:3.9}
\lim_{ A^{\beta}_z \ni x \rightarrow z} \frac{u(x)}{h(x)}
\mbox{ exists for every } \beta >\frac{1-\kappa}{\kappa}.
\end{equation}
\end{thm}

\pf Without loss of generality, we assume  $\nu(\partial D) = 1$. It
 is easy to see that $u(Y^D_t)/h(Y^D_t)$ is a non-negative
supermartingale with respect to $\P^h_{x_0}$. In fact, since $u$ is
non-negative superharmonic for $Y^D$,  $u$ is excessive for $Y^D$.
In particular, $\E_x[  u(Y^D_t)] \le u(x)$ for every x $\in D$. So
by the Markov property for conditioned process, we have for every
$t, s >0$
$$
\E^h_{x_0}\left[  \frac{u(Y^D_{t+s})}{h(Y^D_{t+s})} \,\big|\, {\cal
F}_s\right] \,=\, \E^h_{Y^D_s}\left[
\frac{u(Y^D_{t})}{h(Y^D_{t})}\right]
\,=\,\frac1{h(Y^D_s)}\E_{Y^D_s}\left[  u(Y^D_{t})\right]\, \le\,
\frac{u(Y^D_s)}{h(Y^D_s)}.
$$
Therefore the martingale convergence theorem gives
$$\lim_{t \uparrow \tau^Y_D} \frac{u(Y^D_t)}{h(Y^D_t)}
\mbox { exists and is finite }\P^h_{x_0}\mbox{-a.s. }. $$
Applying
Proposition \ref{prop:3.3}, we have
$$
\int_{\partial D} \P^z_{x_0} \left(\lim_{t \uparrow \tau^Y_D}
\frac{u(Y^D_t)}{h(Y^D_t)} \mbox { exists and is finite}\right)\,
\nu(dz)\,=\, 1.
$$
Thus, for $\nu$-a.e. $z \in \partial D$
\begin{equation}\label{eqn:3.10}
\P^z_{x_0} \left(\lim_{t \uparrow \tau^Y_D}
\frac{u(Y^D_t)}{h(Y^D_t)} \mbox { exists and is finite} \right) = 1.
\end{equation}

We are going to show that (\ref{eqn:3.9})
holds for $z \in \partial D$ satisfying (\ref{eqn:3.10}).
Fix $z \in \partial D$ satisfying (\ref{eqn:3.10})
and fix a $\beta >(1-\kappa)/\kappa$.
Let
$$
l \,:=\, \limsup_{ A^{\beta}_z \ni y \rightarrow z}
\frac{u(y)}{h(y)},
$$
 and assume $l < \infty.$
Then for any $\eps > 0$ there exists a sequence $\{ y_k
\}^{\infty}_{k=1} \subset A_z^{\beta} $ such that $u(y_k)/h(y_k) >
l/(1+\eps )$ and $y_k \rightarrow z.$ (see Lemma 3.9 \cite{K2}).
Without loss of generality, we assume $|y_k-z| < r_0$. Since
$\rho_D(y_k) \le |y_k-z| < r_0$, by Proposition \ref{prop:3.7},
there is $\lambda_0 = \lambda_0 (\eps, \alpha, d)>0 $ such that
\begin{equation}\label{eqn:3.11}
\frac{u(w)}{h(w)} ~\geq~ \frac{u(y_k)}{(
1+\eps )^2h(y_k)}~ >~ \frac{l}{(1+\eps )^3}
\end{equation}
for every $ w \in B^{\lambda_0}_{y_k}=B(y_k, \lambda_0 \rho_D (y))$.

On the other hand,
$$
\P^z_{x_0} \left(T_{B^{\lambda_0}_{y_k}} < \tau^Y_D ~\mbox { i.o.}
\right) \, \geq\, \liminf_{k \rightarrow \infty} \,\P^z_{x_0}
\left(T_{B^{\lambda_0}_{y_k}} < \tau^Y_D \right) \,\geq\, c\,  \,>\,
0 .
$$
But $\{ T_{B^{\lambda_0}_{y_k}} < \tau^Y_D ~
\mbox { i.o.} \}$ is shift-invariant.
Therefore by Proposition \ref{prop:3.5}
\begin{equation}\label{eqn:3.12}
\P^z_{x_0} \left( Y^D_t \mbox{ hits infinitely
many } B^{\lambda_0}_{y_k} \right)
\,=\,\P^z_{x_0} \left(T_{B^{\lambda_0}_{y_k}}
< \tau^Y_D ~~\mbox { i.o.} \right)
\,=\,1 .
\end{equation}
>From (\ref{eqn:3.10})-(\ref{eqn:3.12}), we have
$$
\lim_{t \uparrow \tau_D^Y} \frac{u(Y^D_t)}{h(Y^D_t)} \,\geq\,
\frac{l}{(1+\eps )^3}, \qquad\P^z_{x_0} \mbox{-a.s. for every } \eps
> 0 .$$ Letting $\eps \downarrow 0$,
\begin{equation}\label{eqn:3.13}
\lim_{t \uparrow \tau^Y_D} \frac{u(Y^D_t)}{h(Y^D_t)} \,\geq
\,\limsup_{ A^{\beta}_z \ni y \rightarrow z} \frac{u(y)}{h(y)}
\qquad \P^z_{x_0} \mbox{-a.s. }.
\end{equation}
If $l=\infty$, then for any $M > 1$, there exists a
sequence $\{ y_k \}^{\infty}_{k=1} \subset A_z^{\beta} $
such that $u(y_k)/h(y_k)> 4 M$, $y_k \rightarrow z$
and $\rho_D(y_k) < r_0$.
By Proposition \ref{prop:3.7}, there is $\lambda_1 =
\lambda_1 (M, \alpha, d) >0 $ such that
$$
\frac{u(w)}{h(w)} ~\geq~ \frac{M^2u(y_k)}{(M+1)^2 h(y_k)} \,> \,M
$$
for every $ w \in B^{\lambda_1}_{y_k}$.
So similarly we have
$$
\lim_{t \uparrow \tau^Y_D} \frac{u(Y^D_t)}{h(Y^D_t)}\,
>\, M,  \qquad\P^z_{x_0} \mbox{-a.s. }
$$
for every $M >1$, which is a contradiction because the above limit
is finite $\P^z_{x_0}$-a.s.. Therefore $l < \infty$.

Now let
$$m := \liminf_{ A^{\beta}_z \ni y \rightarrow z}
\frac{u(y)}{h(y)} < \infty.$$
Then for any $\eps > 0$, there exists a sequence $\{ z_k \}^{\infty}_{k=1}
\subset A_z^{\beta} $ such that $u(z_k)/h(z_k)<
m(1+\eps )$, $z_k \rightarrow z$ and $\rho_D(z_k) < r_0$.
By Proposition \ref{prop:3.7},
\begin{equation}\label{eqn:3.11A}
\frac{u(w)}{h(w)} \,\leq\, (1+\eps )^2\,
\frac{u(z_k)} {h(z_k)} \,<\, (1+\eps )^3m
\end{equation}
for every $ w \in B^{\lambda_0}_{z_k}$.
Similarly we have
\begin{equation}\label{eqn:3.12A}
\P^z_{x_0} \left( Y^D_t \mbox{ hits infinitely
many } B^{\lambda_0}_{z_k} \right)
=1 .
\end{equation}
>From (\ref{eqn:3.10}), (\ref{eqn:3.11A}) and
(\ref{eqn:3.12A}), by letting $\eps \downarrow 0$ we have
\begin{equation}\label{eqn:3.14}
\lim_{t \uparrow \tau^Y_D} \frac{u(Y^D_t)}{h(Y^D_t)} \,\leq\,
\liminf_{ A^{\beta}_z \ni y \rightarrow z} \frac{u(y)}{h(y)},
\qquad\P^z_{x_0} \mbox{-a.s. }.
\end{equation}

We conclude from (\ref{eqn:3.13}) and (\ref{eqn:3.14}) that
$$\lim_{ A^{\beta}_z \ni y \rightarrow z}
\frac{u(y)}{h(y)} \mbox{ exists and is finite
for }\nu\mbox{-a.e. } z \in \partial D.$$ \qed

\begin{remark}
{\rm Since constant functions in $\R^d$ are harmonic for $Y$ in $D$,
one can easily see that the above theorem is also
true for every harmonic function $u$ for $Y$
in $D$ either bounded from below or above.
}
\end{remark}

If $u$ and $h$ are harmonic functions for $Y^D$
and $u/h$ is bounded, then $u$ can be recovered
from non-tangential boundary limit values of $u/h$.

\begin{thm}\label{rep}
If $u$ is a harmonic functions for $Y^D$ and $u/h$ is bounded for
some positive harmonic function $h$ for $Y^D$ with
the Martin measure $\nu$,
then for every $x \in D$
$$
u(x)\,=\,h(x)\, \E^h_x\left[\varphi_u \left(\lim_{ t \uparrow
\tau^Y_D} Y^D_t \right)\right],
$$
where
$$
\varphi_u(z):= \lim_{ A^{\beta}_z \ni x \rightarrow z}
\frac{u(x)}{h(x)},\quad\beta >\frac{1-\kappa}{\kappa},
$$
which is well-defined for $\nu$-a.e. $z \in \partial D$.
If we further assume that $u$ is   positive in $D$, then
$$
u(x)= \int_{\partial D} M^Y_D(x,w)\, \varphi_u(w)\, \nu(dw).
$$
That is, $\varphi_u(z)$ is
Radon-Nikodym derivative of the (unique) Martin measure
$\mu_u$ with respect to $\nu$.
\end{thm}

\pf
Without loss of generality, we can assume $u$ is positive and bounded.
 Take an increasing sequence of smooth open sets
$\{D_m\}_{m \ge 1}$ such that $\overline{D_m}
\subset D_{m+1}$ and $\cup^{\infty}_{m=1} D_m = D$.
Let $\tau_m:=\tau^Y_{D_m}$ and $\tau:= \tau^Y_{D}$.
Theorem \ref{thm:son1} implies that
\begin{eqnarray*}
1&=&\P_{x_0}^{z}\left(\lim_{m \rightarrow \infty}
\left(\frac{u}{h}\right)\left(Y^D_{\tau_{m}}\right) =
\lim_{t\uparrow \tau} \left(\frac{u}{h}\right)\left(Y^D_{t}\right)=
\lim_{ A^{\beta}_z \ni x \rightarrow z} \frac{u(x)}{h(x)} \right)\\
&=&\P_{x_0}^{z}\left(\lim_{m \rightarrow \infty}
\left(\frac{u}{h}\right)\left(Y^D_{\tau_{m}}\right) =
\varphi_u(z),~ \lim_{t\uparrow \tau} Y^D_t =z \right)\\
&=&\P_{x_0}^{z}\left(\lim_{m \rightarrow \infty}
\left(\frac{u}{h}\right)\left(Y^D_{\tau_{m}}\right)
 = \varphi_u\left(\lim_{t\uparrow \tau} Y^D_t\right)\right)
\end{eqnarray*}
for $\nu$-a.e. $z \in \partial D$. By Propositions \ref{prop:3.3}
and \ref{prop:3.5}
\begin{equation}\label{df}
\lim_{m \rightarrow \infty}
\left(\frac{u}{h}\right)\left(Y^D_{\tau_{m}}\right) \,=\,
 \varphi_u \left( \lim_{t\uparrow \tau} Y^D_t \right),
\quad\P^h_x \mbox{-a.s. for every } x \in D.
\end{equation}

On the other hand, the harmonicity of $u$ implies that for every $m
\ge1$,
$$ \frac{u(x)}{h(x)} \,=\,\frac{1}{h(x)}\E_x \left[u\left(Y^D_{\tau_{m}}\right) \right]
\,=\, \E^h_x \left[
\left(\frac{u}{h}\right)\left(Y^D_{\tau_{m}}\right) \right].
$$
Therefore, by the bounded convergence theorem and (\ref{df}), we
have
$$
\frac{u(x)}{h(x)}
\,=\,
\lim_{m \rightarrow \infty}
 \E^h_x \left[
\left(\frac{u}{h}\right)\left(Y^D_{\tau_{m}}\right) \right]
 \,=\,\
E^h_x\left[\lim_{m \rightarrow \infty}
\left(\frac{u}{h}\right)\left(Y^D_{\tau_{m}}\right) \right] \,=\,
\E^h_x\left[\varphi_u\left(\lim_{t\uparrow \tau} Y^D_t\right)\right]
$$ for every $x \in D.$ By Theorem \ref{thm:son2},
$$
u(x)= \int_{\partial D} M^Y_D(x,w) \,\varphi_u(w) \,\nu(dw).
$$
\qed

Through an argument similar to the one in Remark 3.19 in \cite{K2},
one can easily check that the above theorem is not true without the
boundedness assumption.

Now suppose that $d=2$, $D=B:=B(0,1)$, $x_0=0$ and
$\sigma_1$ is the normalized surface measure on $\partial B$.
It is showed in \cite{K2} that the Stolz domain is the best
possible one for relative Fatou type theorem for symmetric
stable processes. The proofs there (Lemma 3.22 and Theorem 3.23 in
\cite{K2}) used
Martin Kernel estimates for symmetric
stable process in $B$ and classical argument by Littlewood \cite{L}.
Thus using Martin Kernel estimates for truncated stable
process (Theorem \ref{t:MKE}),
we can show that the Stolz domain is the best
possible one for relative Fatou theorem for truncated
stable processes. We skip the proofs since they are almost
identical to the ones of
Lemma 3.22 and Theorem 3.23 in
\cite{K2}.

\begin{lemma}\label{counter_a}
Suppose that
$$h(x):=\int_{ \partial B } M^Y_B(x, w) \sigma_1(dw)$$
and that $U$ is a measurable function on $\partial B$
such that $0 \le U \le 1.$ Let
$$
u(x):=\int_{ \partial B } M^Y_B(x, w) U(w) \sigma_1(dw) = \frac1{2
\pi} \int^{2 \pi}_0 M^Y_B(x, e^{i \theta}) U( e^{i \theta}) d
\theta,
$$
where $x \in B.$ Suppose that  $0 < \lambda <
\pi $ and $U(e^{i \theta})=1$ for $\theta_0 - \lambda \le \theta \le
\theta_0 + \lambda.$ Then there exists a $\rho = \rho (\eps,
\alpha)$ such that
$$
1-\eps \le \frac{u(\rho e^{i \theta_0})}{h(\rho e^{i \theta_0})} \le
1,  \qquad\mbox{if } \rho > 1- \lambda \rho.
$$
\end{lemma}

A curve $C_0$ is called a tangential curve in $B$
which ends on $\partial B$ if $C_0 \cap \partial
B = \{w_0\} \in \partial B$, $C_0 \setminus \{w_0\}
\subset B$ and there are no $r > 0$ and $ \beta >1$
such that $C_0 \cap B(w_0,r) \subset A^{\beta}_{w_0} \cap B(w_0,r)$.

\begin{thm}\label{counter_thm}
Suppose that
$$
h(x):=\int_{ \partial B } M^Y_B(x, w) \sigma_1(dw).
$$
Let $C_0$ be a tangential curve in $B$ which ends on $\partial B$
and let $C_{\theta}$ be the rotation of $C_{0}$ about $x_0$ through
an angle $\theta.$
Then there exists a positive harmonic function $u$ for $Y$ in
$B:=B(x_0,1)$ such that for a.e. $\theta \in [0, 2\pi]$
with respect to Lebesgue measure,
$$
\lim_{|x| \rightarrow 1, x \in C_{\theta}}
\frac{u(x)}{h(x)} \mbox{ does not exist}.
$$
\end{thm}

\vspace{.3cm} \noindent {\bf Acknowledgment}: We thank
Moritz Kassmann for point out a mistake in a previous version
of this paper.
We also thank an anonymous referee for his many valuable comments
and for suggesting a simpler proof of Theorem \ref{t:ge}.
\vspace{.1in}
\begin{singlespace}

\end{singlespace}
\end{doublespace}

\end{document}